\documentclass[a4paper,10pt]{article}    

 \usepackage{amsmath}
\usepackage{amsfonts}
\usepackage{amssymb}
\usepackage{amsthm}
\usepackage{graphicx}
\usepackage{subfigure}
\usepackage{geometry}
\usepackage[utf8]{inputenc}
\usepackage{float}
\usepackage{pstricks-add}
\usepackage{tikz}
\usetikzlibrary{matrix}     
\usepackage{psfrag}
 
 \newcommand{\ds}{\displaystyle}
 
 \newtheorem{lemma}{Lemma}
\newtheorem{theorem}{Theorem}
\newtheorem{proposition}{Proposition}

\newtheorem{hypothesis}{Hypotheses}

\def\R{\mathbb R}
\def\pa{\partial}
\def\Div{\mathrm{div}}

\def\N{\mathbb N}
\def\R{\mathbb R}
\def\T{\mathcal T}
\def\M{\mathcal M}
\def\P{\mathcal P}
\def\mK{{\rm m}(K)}
\def\msig{{\rm m}(\sigma)}

\def\FKsig{{\mathcal F}_{K,\sigma}}

\def\GKsig{{\mathcal G}_{K,\sigma}}

\def\nKsig {\nu_{K,\sigma}}

\def\ts{\tau_{\sigma}}

\def\E{\mathcal E}
\def\Eint{\E_{int}}
\def\EK{\E_{K}}
\def\EKint{\E_{K,int}}
\def\Eext{\E_{ext}}

\def\EextD{\E_{ext}^D}
\def\EKextD{\E_{K,ext}^D}
\def\EKextN{\E_{K,ext}^N}

\def\somK{\ds\sum_{K\in\T}}
\def\somsigint{\ds\sum_{\substack{\sigma \in \E_{int}\\ \sigma=K|L}}}

\def\somsig{\ds\sum_{\substack{\sigma \in \E\\ (K=K_{\sigma})}}}
\def\somsige{\ds\sum_{\sigma \in \E}}

\def\NKn{N_K^n}
\def\PKn{P_K^n}
\def\NKnp{N_K^{n+1}}
\def\PKnp{P_K^{n+1}}

\def\NKeq{N_K^{eq}}
\def\PKeq{P_K^{eq}}

\def\drNnp{dr(N_{K}^{n+1},N_{K,\sigma}^{n+1})}
\def\drPnp{dr(P_{K}^{n+1},P_{K,\sigma}^{n+1})}

\title{Exponential decay of a finite volume scheme to the thermal equilibrium for drift--diffusion systems}

\author{M. Bessemoulin-Chatard\\ 
Univ. Nantes, CNRS, UMR 6629 - Laboratoire Jean Leray, F-44000 Nantes \and C. Chainais-Hillairet\\
Univ. Lille, CNRS, UMR 8524 - Laboratoire Paul Painlev\'e, F-59000 Lille
}

\begin{document}

\maketitle

\begin{abstract}
In this paper, we study the large--time behavior of a numerical scheme discretizing drift--diffusion systems for semiconductors. The numerical method is finite volume in space, implicit in time, and the numerical fluxes are a generalization of the classical Scharfetter--Gummel scheme which allows to consider both linear or nonlinear pressure laws.\\
We study the convergence of approximate solutions towards an approximation of the thermal equilibrium state as time tends to infinity, and obtain a decay rate by controlling the discrete relative entropy with the entropy production. This result is proved under assumptions of existence and uniform-in-time $L^{\infty}$ estimates for numerical solutions, which are then discussed. We conclude by presenting some numerical illustrations of the stated results. 
\end{abstract}

\section{Introduction}
The Van Roosbroeck's drift--diffusion system is a fundamental model for the mathematical description and numerical simulation of semiconductor devices. It consists of two parabolic convection--diffusion--reaction equations for the carrier densities (electrons and holes), and a Poisson's equation for the electrostatic potential. Global existence and uniqueness results have been obtained for this model under natural assumptions \cite{Gajewski1985,Gajewski1986,Mock1974}. Moreover, this system is shown to be dissipative. Indeed, it admits a Lyapunov functional, which may be physically interpreted as an energy. Furthermore, it has been proved using this energy functional that the solution of the Van Roosbroeck system converges at an exponential rate to the thermal equilibrium state if the boundary conditions are in thermal equilibrium \cite{Gajewski1996,Gajewski1989}.

The classical drift--diffusion model is based on Boltzmann statistics. More precisely, it means that the statistical distribution function describing the dependence of the carrier densities on the chemical potentials is the exponential function. However, this choice may fail to describe relevantly the physical reality in some cases (for example in case of high carrier densities).Then other statistics have to be considered, like Fermi-Dirac statistics for instance \cite{Gajewski1989}. This leads to a modification of the diffusive terms, which become nonlinear. Existence and uniqueness of weak solutions to a nonlinear drift-diffusion model have been proved in \cite{Juengel1994}. Dissipativity and convergence to the thermal equilibrium for large time have also been established for this generalized model \cite{Juengel1995}. Let us underline that more recently, driven by applications like organic semiconductors, there is an increased interest in drift--diffusion models with arbitrary statistical distribution functions \cite{Foster2014,Mensfoort2008}. 

From a numerical point of view it is essential to consider numerical schemes which preserve the main qualitative properties of the continuous system, such as positivity of the densities, dissipativity and consistency with thermal equilibrium. In the case of Boltzmann statistics, the Scharfetter--Gummel scheme \cite{IlIn1969,Scharfetter1969} is widely used. It exploits the exponential dependence on the chemical potential and allows to recover the correct large time behavior. Various extensions of the Scharfetter--Gummel scheme have been suggested to account the diffusion enhancement induced by non Boltzmann statistics \cite{Juengel1995a,Purbo1989,Stodtmann2012}. Unfortunately, they are not thermodynamically consistent. More recently, a consistent generalization in the spirit of the original Scharfetter--Gummel scheme was proposed \cite{Eymard2006} and applied to simulate organic semiconductors' behavior \cite{Koprucki2013,Koprucki2013a}. This method leads to solve a nonlinear boundary value problem at each interface.

In this paper, we preferentially focus on another extension of the Scharfetter--Gummel scheme using a proper average of the nonlinear diffusion \cite{Bessemoulin-Chatard2012} which guarantees thermodynamic consistency. An alternative interpretation of this scheme is given in \cite{Koprucki2015} and applied to a very general class of statistical distribution functions arising in organic semiconductors modeling.  
Here our aim is to study the large--time behavior of an implicit in time and finite volume in space discretization of the drift--diffusion system, with a Scharfetter--Gummel approximation of the convection--diffusion fluxes. Our proof is based on the entropy--dissipation method \cite{Arnold2004}. This point of view was already adopted in several articles \cite{Chainais-Hillairet2007,Gajewski1996,Glitzky2008,Glitzky2011,Glitzky2009}. The crucial point to obtain an exponential decay rate of the approximate solutions towards the equilibrium is the control of the relative entropy by the entropy production. The choice of Scharfetter--Gummel type fluxes for the discretization of the convection--diffusion fluxes is essential at this step.


\subsection{The drift--diffusion system and the thermal equilibrium}

Let $\Omega$ be an open bounded subset of $\mathbb{R}^{d}$ ($d\geq 1$) corresponding to the geometry of a semiconductor device and $T>0$. This device can be described by the so-called drift--diffusion system. This system consists of two continuity equations for the electron density $N$ and the hole density $P$, and a Poisson equation for the electrostatic potential $\Psi$. It writes for all $(x,t)\in\Omega\times [0,T]$:
\begin{gather}
\pa_{t}N+\Div(\mu_{N}(-\nabla r(N)+N\nabla\Psi))=-R(N,P),\label{systemDD-N}\\
\pa_{t}P+\Div(\mu_{P}(-\nabla r(P)-P\nabla\Psi))=-R(N,P),\label{systemDD-P}\\
-\lambda^{2}\Delta\Psi=P-N+C.\label{systemDD-Psi}
\end{gather}
The given function $C(x)$ is the doping profile describing fixed background charges. The dimensionless physical parameters $\mu_{N}$, $\mu_{P}$ and $\lambda$ are the rescaled mobilities of electrons and holes, and the rescaled Debye length respectively. The definition of $r$ depends on the statistics chosen to describe the relation between the densities and the chemical potentials. The usual considerations on which the isentropic hydrodynamic model are based suggest a pressure of the form:
$$ r(s)=s^{\alpha},\quad\alpha\geq 1.$$
The linear case, where $\alpha=1$, is the isothermal model, corresponding to Boltzmann statistics. 
The system \eqref{systemDD-N}--\eqref{systemDD-Psi} is supplemented with initial conditions:
\begin{equation}\label{CI}
N(x,0)=N_{0}(x), \quad P(x,0)=P_{0}(x), \quad x\in\Omega,
\end{equation}
and with mixed boundary conditions: Dirichlet boundary conditions on the ohmic contacts and homogeneous Neumann boundary conditions on the insulated boundary segments. More precisely, the boundary $\pa \Omega$ is split into $\pa \Omega=\Gamma^{D}\cup\Gamma^{N}$ with $\Gamma^{D}\cap\Gamma^{N}=\emptyset$, and the boundary conditions write:
\begin{equation}\label{BCD}
N(\gamma,t)=N^{D}(\gamma), \, P(\gamma,t)=P^{D}(\gamma), \, \Psi(\gamma,t)=\Psi^{D}(\gamma), \,\, (\gamma,t)\in\Gamma^{D}\times [0,T],
\end{equation}
\begin{equation}\label{BCN}
(\nabla r(N)\cdot \nu)(\gamma,t)=(\nabla r(P)\cdot \nu)(\gamma,t)=(\nabla \Psi\cdot \nu)(\gamma,t)=0,\,\, (\gamma,t)\in\Gamma^{N}\times [0,T],
\end{equation}
where $\nu$ is the unit normal to $\pa\Omega$ outward to $\Omega$. In this paper, we assume that the mobilities are constant and equal: $\mu_{N}=\mu_{P}=1$, and we need the following general assumptions:
\begin{hypothesis}\label{HYP}
The domain $\Omega$ is an open bounded polygonal (or polyhedral) subset of $\mathbb{R}^d$ ($d\geq 1$) and $\pa \Omega=\Gamma^D \cup\Gamma^N$, with $\Gamma^D\cap\Gamma^N=\emptyset$ and $\text{m}(\Gamma^D)>0$.
The doping profile $C$ belongs to $L^\infty(\Omega)$.
The boundary conditions $N^D$, $P^D$ and $\Psi^D$ are traces of some functions defined on the whole domain $\Omega$, still denoted by $N^D$, $P^D$ and $\Psi^D$.
Furthermore, we assume that
\begin{gather}
N_{0},\,P_{0}\in L^\infty(\Omega),\label{HYP-CI}\\
N^D,\,P^D\in L^\infty\cap H^1(\Omega),\,\Psi^D\in H^1(\Omega),\label{HYP-BC}\\
\exists M>0,\,m\geq 0 \text{ such that } m\leq N_{0},P_{0},N^D,P^D\leq M\,\text{a.e. on }\Omega.\label{HYP-M-CI-BC}
\end{gather}
\end{hypothesis}

For the drift--diffusion model with linear pressure $r=Id$, the recombination--generation rate can usually be written under the following form \cite{Markowich1990}:
\begin{equation}\label{forme-R}
R(N,P)=R_{0}(N,P)(NP-1).
\end{equation}
This general form includes in particular the Shockley--Read--Hall term:
$$R_{SRH}(N,P)=\frac{NP-1}{\tau_{P}N+\tau_{N}P+\tau_{C}}, \quad \tau_{P},\,\tau_{N},\,\tau_{C}>0,$$
or the Auger recombination:
$$R_{AU}=(C_{N}N+C_{P}P)(NP-1).$$
In the case of a nonlinear pressure, the Shockley--Read--Hall term cannot be taken anymore. Some recent works about organic electronic devices (see \cite{Gruber2011} for instance) include several propositions for the modeling of generation and recombination processes in the nonlinear case. However, this leads to highly nonlinear and intricate source terms, and we choose to consider either the linear case with $R$ defined by \eqref{forme-R}, or the nonlinear case with $R=0$. More precisely, we will need to assume that either
\begin{hypothesis}\label{HYP-LIN}
\begin{gather}
  r=Id, \label{HYP-LIN-1}\\
R(N,P)=R_{0}(N,P)(NP-1), \text{ with } R_{0} \text{ continuous and nonnegative}, \label{HYP-LIN-2}\\
 N^DP^D=1,\label{HYP-LIN-3}
\end{gather}
\end{hypothesis}
\noindent or
\begin{hypothesis}\label{HYP-NON-LIN}
\begin{gather}
 r\in\mathcal{C}^{1}(\mathbb{R}), \quad r(0)=r'(0)=0, \quad r'(s)\geq c_{0}s^{\alpha-1}, \quad \alpha>1, \label{HYP-NON-LIN-1}\\
 R=0.  \label{HYP-NON-LIN-2}
\end{gather}
\end{hypothesis}

Existence and uniqueness of weak solutions to the drift--diffusion system have been studied in \cite{Gajewski1985,Gajewski1986,Mock1974} for the isothermal model, whereas the nonlinear case is considered in \cite{Juengel1994}.
The large time behavior of the isothermal drift--diffusion system \eqref{systemDD-N}--\eqref{BCN} has been studied in~\cite{Gajewski1996}. It has been proven that the solution to the transient system converges to the thermal equilibrium state as $t\rightarrow +\infty$ if the boundary conditions \eqref{BCD} are in thermal equilibrium. This result was extended to the degenerate case with nonlinear diffusivities in \cite{Juengel1995}.\\
More precisely, the thermal equilibrium is a particular steady--state for which electron and hole currents vanish, namely
$$-\nabla r(N)+N\nabla\Psi=-\nabla r(P)-P\nabla\Psi=0.$$
The existence of a thermal equilibrium has been studied in the case of a linear pressure in \cite{Markowich1990,Markowich1986} and in the nonlinear case in \cite{Markowich1993}. Let us introduce the enthalpy function $h$ defined by:
$$h(s)=\int_{1}^{s}\frac{r'(\tau)}{\tau}d\tau,$$
and the generalized inverse $g$ of $h$, defined by:
\begin{equation*}
g(s)=\left\{\begin{array}{lcl}
h^{-1}(s) & \text{ if } & h(0^{+})<s<\infty,\\
0 & \text{ if } & s\leq h(0^{+}),
\end{array}\right.
\end{equation*}
where we have implicitly assumed that $h(+\infty)=+\infty$. In the isothermal case, we simply have $h=\log$ and $g=\exp$. \\
If the Dirichlet boundary conditions satisfy $N^{D}$, $P^{D}>0$ and
\begin{equation}\label{hyp-compatibility-BC}
h(N^{D})-\Psi^{D}=\alpha_{N} \text{ and } h(P^{D})+\Psi^{D}=\alpha_{P} \text{ on }\Gamma^{D},
\end{equation}
the thermal equilibrium is defined for $x\in\Omega$ by:
\begin{gather}
-\lambda^{2}\Delta\Psi^{eq}=g(\alpha_{P}-\Psi^{eq})-g(\alpha_{N}+\Psi^{eq})+C,\label{eqtherm-Psi}\\
N^{eq}=g(\alpha_{N}+\Psi^{eq}),\label{eqtherm-N}\\
P^{eq}=g(\alpha_{P}-\Psi^{eq}),\label{eqtherm-P}
\end{gather}
with the boundary conditions \eqref{BCD}--\eqref{BCN}.\\
As mentioned in Hypotheses \ref{HYP-LIN}, we furthermore assume in the linear case that the Dirichlet boundary conditions satisfy the mass action law \eqref{HYP-LIN-3}: $N^D P^D=1$. At the thermal equilibrium, $R(N^{eq},P^{eq})=0$ must hold, which implies in view of the form of $R$ \eqref{HYP-LIN-2} that $N^{eq}P^{eq}=1$, and finally that $\alpha_{N}+\alpha_{P}=0$.\\
The proof of convergence to the thermal equilibrium is based on the entropy method, described for instance in the review paper \cite{Arnold2004}. This method consists of looking for a nonnegative Lyapunov functional, called entropy, and its nonnegative production, connected within an entropy--entropy production estimate. It provides the convergence in relative entropy of the evolutive solution towards the equilibrium state. Moreover, if the relative entropy is controlled with the entropy production, one can compute a convergence rate. This method has been widely applied to many different systems ; see for instance \cite{Arnold2001} for Fokker--Planck type equations, or \cite{Carrillo2001} for degenerate parabolic problems. \\
Here the relative entropy functional is the deviation of the total energy (sum of the internal energies for the electron and hole densities and the energy due to the electrostatic potential) from the thermal equilibrium:
\begin{multline}\label{defEcontinu}
\mathbb{E}(t)=\int_{\Omega}\left(\vphantom{\frac{1}{2}}H(N(t))-H(N^{eq})-h(N^{eq})(N(t)-N^{eq})\right.\\
+ H(P(t))-H(P^{eq})-h(P^{eq})(P(t)-P^{eq})\\
\left.+\frac{\lambda^{2}}{2}|\nabla(\Psi(t)-\Psi^{eq})|^{2}\right)dx,
\end{multline}
with $H(x)=\int_{1}^{x}h(s)\,ds$, and the entropy production functional is given by
\begin{multline}\label{defIcontinu}
\mathbb{I}(t)=\int_{\Omega}\left(\vphantom{\frac{1}{2}}N|\nabla(h(N)-\Psi)|^{2}+P|\nabla(h(P)+\Psi)|^{2}\right)dx\\
 +\int_{\Omega}R(N,P)\left(h(N)+h(P)-h(N^{eq})-h(P^{eq})\right)dx.
\end{multline}
In the nonlinear case, we assume that $R=0$, and then the last term of $\mathbb{I}(t)$ vanishes, whereas in the linear case with recombination--generation rate of the form \eqref{HYP-LIN-2}, we obtain that
\begin{gather*}
\mathbb{I}(t)=\int_{\Omega}\left(\vphantom{\frac{1}{2}}N|\nabla(\log(N)-\Psi)|^{2}+P|\nabla(\log(P)+\Psi)|^{2}\right)dx\\
+\int_{\Omega}R_{0}(N,P)(NP-1)\log(NP),
\end{gather*}
where the last term is clearly nonnegative.

The entropy--entropy production inequality writes:
\begin{equation}\label{ineg_EI_continu_1}
0 \leq \mathbb{E}(t)+\int_{0}^{t}\mathbb{I}(s)\,ds\leq\mathbb{E}(0).
\end{equation}
Exponential decay towards the thermal equilibrium for the drift--diffusion model has been proved in \cite{Gajewski1986} for the Boltzmann statistics, and extended to the Fermi-Dirac statistics in \cite{Gajewski1989}. In \cite{Gajewski1996}, the long--time behavior of the model with magnetic field is studied. Large--time behavior of reaction--diffusion systems for a finite number of charged species has been investigated in \cite{Gajewski1996a,Glitzky1996,Glitzky2009}. Finally, the convergence towards the thermal equilibrium state for the drift--diffusion system in the nonlinear degenerate case is proved in \cite{Juengel1995}, but without any rate.


\subsection{Outline of the paper}

The outline of the paper is as follows. In Section \ref{sec-discret}, we introduce the discrete framework. It includes the description of the considered numerical schemes as well as the definition of the discrete relative entropy and the corresponding entropy production. We conclude this section by establishing a technical property of the generalized Scharfetter--Gummel fluxes in Lemma \ref{lem-flux}. Section \ref{sec-exp-decay} is devoted to the detailed proof of our main result, which is stated in Theorem \ref{thrm:decayE}. The decay of numerical solutions towards the equilibrium is studied either under Hypotheses~\ref{HYP-LIN} corresponding to the isothermal case with nonzero recombination--generation rate, or under Hypotheses~\ref{HYP-NON-LIN} corresponding to a rather general nonlinear pressure, without recombination--generation rate. We also assume existence and uniform $L^{\infty}$ estimates for numerical solutions to establish the large--time behavior of the scheme. These assumptions are then discussed in Section \ref{sec-existence}, where we distinguish the nonlinear case from the isothermal case. It appears that the uniform-in-time $L^{\infty}$ estimates needed to prove Theorem \ref{thrm:decayE} are only obtained in the case of a zero doping profile. However in the last section, we present some numerical results and observe an exponential convergence towards a  steady-state even when this condition is not satisfied.


\section{Presentation of the discrete setting}
\label{sec-discret}


\subsection{Definition of the numerical schemes}

In this subsection, we present the finite volume schemes for the time evolution drift--diffusion system \eqref{systemDD-N}--\eqref{BCN} and for the thermal equilibrium \eqref{eqtherm-Psi}--\eqref{eqtherm-P}. The mesh $\M=(\T,\E,\P)$ of the domain $\Omega$ is given by a family $\T$ of open polygonal (or polyhedral in 3-D) control volumes, a family $\E$ of edges (or faces), and a family $\P=(x_{K})_{K\in\T}$ of points. As it is classical in the finite volume discretization of diffusive terms with two-points flux approximations, we assume that the mesh is admissible in the sense of \cite[Definition 9.1]{Eymard2000}. It implies that the straight line between two neighboring centers of cells $(x_{K},x_{L})$ is orthogonal to the edge $\sigma=K|L$.\\
In the set of edges $\E$, we distinguish the interior edges $\sigma=K|L\in \Eint$ and the boundary edges $\sigma\in\Eext$. Within the exterior edges, we distinguish the Dirichlet boundary edges included in $\Gamma^{D}$ from the Neumann boundary edges included in $\Gamma^{N}$: $\Eext=\Eext^{D}\cup\Eext^{N}$. For a control volume $K\in\T$, we define $\EK$ the set of its edges, which is also split into $\EK=\EKint\cup\EKextD\cup\EKextN$. For each edge $\sigma\in\E$, there exists at least one cell $K\in\T$ such that $\sigma\in\EK$, which will be denoted $K_{\sigma}$. In the case where $\sigma=K|L\in\Eint$, $K_{\sigma}$ can be either equal to $K$ or $L$.\\
For all $\sigma\in\E$, we define $d_{\sigma}=d(x_{K},x_{L})$ if $\sigma=K|L\in\Eint$ and $d_{\sigma}=d(x_{K},\sigma)$ if $\sigma\in\Eext$, with $\sigma\in\EK$. Then the transmissibility coefficient is defined by $\ts=\msig/d_{\sigma}$, for all $\sigma\in\E$.\\
We assume that the mesh satisfies the following regularity constraint:
\begin{equation}\label{regmesh}
\exists\xi >0 \text{ such that }d(x_{K},\sigma)\geq \xi\,d_{\sigma},\quad\forall K\in\T,\quad \forall\sigma\in\EK.
\end{equation}
Let $\Delta t>0$ be the time step. We assume that there exists $\Delta t_{max}>0$ such that
\begin{equation}\label{hyp-dt}
0\leq \Delta t\leq \Delta t_{max}.
\end{equation}
We set $N_{T}=E(T/\Delta t)$ and $t^{n}=n\Delta t$ for all $0\leq n\leq N_{T}$. The size of the mesh is defined by $\text{size}(\T)=\max_{K\in\T}\text{diam}(K)$, and we denote by $\delta=\max(\Delta t,\text{size}(\T))$ the size of the space--time discretization.\\
A finite volume scheme for a conservation law with unknown $u$ provides a vector $u_{\T}=(u_{K})_{K\in\T}\in\R^{\theta}$ (with $\theta=\text{Card}(\T)$) of approximate values and the associate piecewise constant function, still denoted $u_{\T}$:
\begin{equation*}
u_{\T}=\sum_{K\in\T}u_{K}\mathbf{1}_{K},
\end{equation*}
where $\mathbf{1}_{K}$ denotes the characteristic function of the cell $K$. However, since there are Dirichlet conditions on a part of the boundary, we also need to define approximate values for $u$ at the corresponding boundary edges: $u_{\E^{D}}=(u_{\sigma})_{\sigma\in\EextD}\in\R^{\theta^{D}}$ (with $\theta^{D}=\text{Card}(\EextD)$). Therefore, the vector containing the approximate values both in the control volumes and at the Dirichlet boundary edges is denoted by $u_{\M}=(u_{\T},u_{\E^{D}})$. We denote by $X(\M)$ the set of the discrete functions $u_\M$ and by $X_0(\M)$ the subset of $X(\M)$ of the functions vanishing at the boundary $\Gamma^D$:  $u_{\M}=(u_{\T},0_{\E^{D}})$.

For any vector $u_{\M}=(u_{\T},u_{\E^{D}})$, we define for all $K\in\T$ and all $\sigma\in\EK$
\begin{equation*}
u_{K,\sigma}=\left\{\begin{array}{ll}
u_{L}& \text{ if } \sigma=K|L\in\EKint,\\
u_{\sigma}& \text{ if } \sigma\in\EKextD,\\
u_{K} & \text{ if }\sigma\in\EKextN,
\end{array}\right.
\end{equation*}
and 
\begin{equation*}
Du_{K,\sigma}=u_{K,\sigma}-u_{K}, \quad  \quad D_{\sigma}u=|Du_{K,\sigma}|.
\end{equation*}
We can now define the discrete $L^2$ norm and the discrete $H^{1}$-seminorm $|\cdot|_{1,\M}$ on $X(\M)$ by
$$
\begin{aligned}
|u_{\M}|_{1,\M}^{2}&=\somsige\ts(D_{\sigma}u)^{2},\quad\forall u_{\M}\in X(\M),\\
\Vert u_{\M}\Vert_{0,\M}^{2}&=\Vert u_\T\Vert_0^2=\ds\sum_{K\in\T} \mK (u_K)^2 ,\quad\forall u_{\M}\in X(\M),
\end{aligned}
$$
where $\|\cdot\|_{0}$ is the classical $L^2$ norm for piecewise constant functions.

Let us recall in Proposition \ref{prop_Poincare} the discrete counterpart of the Poincar\'e inequality for piecewise constant functions. We refer to \cite[Theorem 4.3]{Bessemoulin-Chatard2015} for a proof of this result in the case of mixed boundary conditions.
\begin{proposition}\label{prop_Poincare}
Let $\Omega$ be an open bounded polyhedral domain of $\R^d$ ($d\geq 1$) and $\pa \Omega=\Gamma^D\cup\Gamma^N$ with $\Gamma^D\cap\Gamma^N=\emptyset$ and $m(\Gamma^D)>0$. Let $\mathcal M=(\T,\E,\P)$ an admissible finite volume mesh of $\Omega$ which satisfies \eqref{regmesh}. There exists a constant $C_P$ depending only on $\Omega$ such that 
\begin{equation}\label{ineq_Poincare}
\Vert u_\M\Vert_{0,\M}\leq \frac{C_P}{\xi^{1/2}} |u_{\M}|_{1,\M} \quad \forall u_\M\in X_0(\M).
\end{equation}
\end{proposition}

\paragraph{The scheme for the transient model.}
We have to define at each time step $0\leq n\leq N_{T}$ the approximate solution $u_{\T}^{n}=(u_{K}^{n})_{K\in\T}$ for $u=N,\,P,\,\Psi$ and the approximate values at the boundary $u_{\E^{D}}^{n}=(u_{\sigma}^{n})_{\sigma\in\Eext^{D}}$ (which in fact does not depend on $n$ since the boundary data do not depend on time).  First of all, we discretize the initial and boundary conditions: 
\begin{gather}
\left(N_{K}^{0},P_{K}^{0}\right)=\frac{1}{\mK}\int_{K}\left(N_{0}(x),P_{0}(x)\right)\,dx,\quad \forall K \in\T, \label{CIdiscret}\\
\left(N_{\sigma}^{D},P_{\sigma}^{D},\Psi_{\sigma}^{D}\right)=\frac{1}{\msig}\int_{\sigma}\left(N^{D}(\gamma),P^{D}(\gamma),\Psi^{D}(\gamma)\right)d\gamma,\quad\forall\sigma\in\EextD, \label{BCdiscret}
\end{gather}
and we define
\begin{equation}\label{BCdiscretn}
N_{\sigma}^{n}=N_{\sigma}^{D},\quad P_{\sigma}^{n}=P_{\sigma}^{D},\quad \Psi_{\sigma}^{n}=\Psi_{\sigma}^{D},\quad \forall\sigma\in\EextD,\quad\forall n\geq 0.
\end{equation}
 Then, we consider a backward Euler in time and finite volume in space discretization of the drift--diffusion system \eqref{systemDD-N}--\eqref{systemDD-Psi}. The scheme writes:\\ $\quad \forall K\in{\mathcal T},\, \forall n\geq 0,$
\begin{align}
&\mK \frac{N_K^{n+1}-N_K^n}{\Delta t}+\ds\sum_{\sigma\in {\mathcal E}_K}{\mathcal F}_{K,\sigma}^{n+1}=-\mK\,R(N_{K}^{n+1},P_{K}^{n+1}),\label{scheme-N}\\
&\mK \frac{P_K^{n+1}-P_K^n}{\Delta t}+\ds\sum_{\sigma\in {\mathcal E}_K}{\mathcal G}_{K,\sigma}^{n+1}=-\mK\,R(N_{K}^{n+1},P_{K}^{n+1}),\label{scheme-P}\\
&-\lambda^2\sum_{\sigma\in {\mathcal E}_K}\tau_{\sigma} D\Psi_{K,\sigma}^{n}=\mK (P_K^{n}-N_K^{n}+C_{K}).\label{scheme-Psi}
\end{align}

It remains to define the numerical fluxes $\FKsig^{n+1}$ and $\GKsig^{n+1}$ which are approximations of $\ds\int_{\sigma} (-\nabla r(N)+N\nabla\Psi)\cdot \nKsig$ and $\ds\int_{\sigma} (-\nabla r(P)-P\nabla\Psi)\cdot \nKsig$ on the interval $[t^n,t^{n+1})$. We choose to discretize simultaneously the diffusive part and the convective part of the fluxes. In the linear case $r=Id$, we use the classical Scharfetter--Gummel fluxes. For all $K\in\T$, for all $\sigma\in\E_K$, we set: 
\begin{align}
&\label{FLUX-SG-N}
{\mathcal F}_{K,\sigma}^{n+1}=
\tau_{\sigma}\left[ B\left(-D\Psi_{K,\sigma}^{n+1}\right)N_K^{n+1}-B\left(D\Psi_{K,\sigma}^{n+1}\right)N_{K,\sigma}^{n+1} \right],\\
&\label{FLUX-SG-P}
{\mathcal G}_{K,\sigma}^{n+1}=
\tau_{\sigma}\left[ B\left(D\Psi_{K,\sigma}^{n+1}\right)P_K^{n+1}-B\left(-D\Psi_{K,\sigma}^{n+1}\right)P_{K,\sigma}^{n+1} \right],
\end{align}
where $B$ is the Bernoulli function defined by 
\begin{equation}\label{Bern}
B(0)=1 \mbox{ and } B(x)=\ds\frac{x}{\exp (x) -1}\quad  \forall x\neq 0.
\end{equation}
These fluxes, introduced in \cite{IlIn1969,Scharfetter1969}, are widely used to approximate the drift--diffusion system in the linear case. They are second order accurate in space \cite{Lazarov1996} and they preserve steady--states. In \cite{Gajewski1996}, the dissipativity of the implicit Scharfetter--Gummel scheme is proved. A proof of the exponential decay of the free energy to its equilibrium value is given in \cite{Glitzky2009} for an implicit time discretization of electro--reaction--diffusion problems, and this result is extended to a fully discrete problem in \cite{Glitzky2008,Glitzky2011}. In \cite{Gaertner2009}, some bounds for discrete steady--states solutions obtained with the Scharfetter--Gummel scheme are established. Moreover, a discrete analog of the entropy--entropy production inequality \eqref{ineg_EI_continu_1} is proved in \cite{Chatard2011}, yielding the long--time behavior of the Scharfetter--Gummel scheme for the linear drift--diffusion system.

In the case of a nonlinear pressure function $r$ satisfying \eqref{HYP-NON-LIN-1},  we use a generalization of the Scharfetter--Gummel fluxes defined in \cite{Bessemoulin-Chatard2012}. For all $K\in\T$, for all $\sigma\in\E_K$, we set: 
\begin{align}
\label{FLUX-SG-N-ext}
{\mathcal F}_{K,\sigma}^{n+1}=
\tau_{\sigma}\,\drNnp & \left[ B\left(\frac{-D\Psi_{K,\sigma}^{n+1}}{\drNnp}\right)N_K^{n+1}\right. \\
&\left. -B\left(\frac{D\Psi_{K,\sigma}^{n+1}}{\drNnp}\right)N_{K,\sigma}^{n+1} \right],\nonumber\\
\label{FLUX-SG-P-ext}
{\mathcal G}_{K,\sigma}^{n+1}=
\tau_{\sigma}\,\drPnp & \left[ B\left(\frac{D\Psi_{K,\sigma}^{n+1}}{\drPnp}\right)P_K^{n+1} \right. \\
&\left.-B\left(\frac{-D\Psi_{K,\sigma}^{n+1}}{\drPnp}\right)P_{K,\sigma}^{n+1} \right],\nonumber
\end{align}
where
\begin{equation}\label{def-dr}
dr(a,b)=\left\{\begin{array}{lcl}
\ds{\frac{h(b)-h(a)}{\log(b)-\log(a)}} & \text{ if } & a, \, b>0, \quad a \neq b,\\
\ds{r'\left(\frac{a+b}{2}\right)} & \text{ elsewhere.}
\end{array}\right.
\end{equation}
This particular choice of $dr$ ensures the preservation of the thermal equilibrium at the discrete level, which is crucial to have a good long time behavior. We notice that in the isothermal case $r=Id$, we recover exactly the classical Scharfetter--Gummel fluxes \eqref{FLUX-SG-N} and \eqref{FLUX-SG-P}. These fluxes were also generalized to a larger class of statistical distribution functions arising in organic semiconductors modeling \cite{Koprucki2015}.

Other extensions of the Scharfetter--Gummel scheme have been proposed. A scheme valid in the case where both convective and diffusive terms are nonlinear is studied in \cite{Eymard2006} and applied to organic semiconductors models in \cite{Koprucki2013,Koprucki2013a}, but this method leads to solve a nonlinear elliptic problem at each interface. We also mention \cite{Juengel1995a,Juengel1997} where fluxes \eqref{FLUX-SG-N-ext}--\eqref{FLUX-SG-P-ext} are considered, but with another definition of $dr$ which does not allow to preserve the thermal equilibrium at the discrete level. A finite volume scheme preserving the long--time behavior of the solutions of the nonlinear drift--diffusion model is introduced in \cite{Chainais-Hillairet2007}, and the convergence of this scheme towards the equilibrium state is proved, based on the control of the discrete energy production.

\paragraph{The scheme for the thermal equilibrium.} We compute an approximation of the thermal equilibrium  $(N_{\T}^{eq},P_{\T}^{eq},\Psi_{\T}^{eq})$ defined by \eqref{eqtherm-Psi}--\eqref{eqtherm-P} with the following finite volume scheme: $\forall K\in{\mathcal T},$
\begin{gather}
-\lambda^2\sum_{\sigma\in {\mathcal E}_K}\tau_{\sigma} D\Psi_{K,\sigma}^{eq}=\mK \left(g(\alpha_{P}-\Psi_{K}^{eq})-g(\alpha_{N}+\Psi_{K}^{eq})+C_{K}\right), \label{scheme-Psi-eq}\\
N_{K}^{eq}=g(\alpha_{N}+\Psi_{K}^{eq}),\label{scheme-N-eq}\\
P_{K}^{eq}=g(\alpha_{P}-\Psi_{K}^{eq}).\label{scheme-P-eq}
\end{gather}
Existence and uniqueness of solution to this nonlinear scheme is studied in \cite{Chainais-Hillairet2007}. 
$L^{\infty}$ estimates on $\Psi_{\T}^{eq}$ are also established in \cite{Chainais-Hillairet2007}. Due to \eqref{scheme-N-eq} and \eqref{scheme-P-eq}, they imply lower and upper bounds on $N_{\T}^{eq}$ and $P_{\T}^{eq}$. The lower bound is positive in the linear case and in the case of  Fermi-Dirac statistics (see \cite{Gajewski1989}) but it may vanish in the general nonlinear case. We assume in  what follows that this lower bound is positive.

\paragraph{The discrete entropy and entropy production.} For $n\in\N$, the discrete relative entropy functional is defined by:
\begin{multline}\label{defEdiscret}
\mathbb{E}^{n}=\somK\mK\left(H(\NKn)-H(\NKeq)-h(\NKeq)(\NKn-\NKeq)\right. \\
\left. + H(\PKn)-H(\PKeq)-h(\PKeq)(\PKn-\PKeq)\right)\\
+\frac{\lambda^{2}}{2}\left|\Psi_{\M}^{n}-\Psi_{\M}^{eq}\right|^{2}_{1,\M}
\end{multline}
and the discrete entropy production functional is defined by:
\begin{multline}\label{defIdiscret}
\mathbb{I}^{n}=\somsig \ts \left[\min(\NKn,N_{K,\sigma}^{n})\left(D_{\sigma}(h(N^{n})-\Psi^{n})\right)^{2}\right.\\ 
+\left.\min(\PKn,P_{K,\sigma}^{n})\left(D_{\sigma}(h(P^{n})+\Psi^{n})\right)^{2}\right]\\
+\somK\mK R(\NKn,\PKn)\left[h(\NKn)+h(\PKn)-h(\NKeq)-h(\PKeq)\right],
\end{multline}
where $\ds \somsig$ means a sum over all the edges $\sigma\in\E$ and $K$ inside the sum is replaced by $K_{\sigma}$ (therefore $\sigma$ is an edge of the cell $K=K_{\sigma}$).

As already mentioned in the continuous framework, the last term in the definition of $\mathbb{I}^n$ is reduced to zero in the nonlinear case where $R=0$ and is nonnegative in the linear case with a recombination--generation rate under the form \eqref{HYP-LIN-2}.


\subsection{Properties of the numerical fluxes}

In Section \ref{sec-exp-decay}, we study the large time behavior of the approximate solution by adapting the entropy--dissipation method to the discrete level. We start by establishing a lemma which will be useful to prove the discrete counterpart of \eqref{ineg_EI_continu_1}. It is a generalization to the nonlinear case of Corollary A.2 in  \cite{Bessemoulin-Chatard2014a}.

\begin{lemma}\label{lem-flux}
For all $K\in\T$ and all $\sigma\in\E_{K}$, the fluxes $\FKsig^{n+1}$ and $\GKsig^{n+1}$ defined by \eqref{FLUX-SG-N-ext} and \eqref{FLUX-SG-P-ext} satisfy:
\begin{align}
 \FKsig^{n+1}\,D\left(h(N^{n+1})- \right. & \left. \Psi^{n+1}\right)_{K,\sigma} \leq  \label{prop-FLUX-N}\\
& -\ts\,\min(\NKnp,N_{K,\sigma}^{n+1})\left(D_{\sigma}\left(h(N^{n+1})-\Psi^{n+1}\right)\right)^2,  \nonumber\\
 \GKsig^{n+1}\,D\left(h(P^{n+1})+\right. & \left.\Psi^{n+1}\right)_{K,\sigma}  \leq  \label{prop-FLUX-P}\\
 & -\ts\,\min(\PKnp,P_{K,\sigma}^{n+1})\left(D_{\sigma}\left(h(P^{n+1})+\Psi^{n+1}\right)\right)^2. \nonumber
\end{align}
\end{lemma}

\begin{proof}
We focus on the proof of \eqref{prop-FLUX-N} because \eqref{prop-FLUX-P} will be proved in the same way, replacing $N$ by $P$ and $\Psi$ by $-\Psi$. Let us define 
\begin{align*}
\mathcal{R}_{K,\sigma}^{n+1} &=
\FKsig^{n+1}\,D\left(h(N^{n+1})-\Psi^{n+1}\right)_{K,\sigma}\\
&\,\,\,+ \ts\,\min(\NKnp,N_{K,\sigma}^{n+1})\left(D_{\sigma}\left(h(N^{n+1})-\Psi^{n+1}\right)\right)^2\\
&=D\left(h(N^{n+1})-\Psi^{n+1}\right)_{K,\sigma} \\
&\quad \times\left(\FKsig^{n+1}+\ts\,\min(\NKnp,N_{K,\sigma}^{n+1})D\left(h(N^{n+1})-\Psi^{n+1}\right)_{K,\sigma}\right).
\end{align*}
Our aim is to prove that $\mathcal{R}_{K,\sigma}^{n+1}\leq 0$. We may write:
\begin{gather*}
 D\left(h(N^{n+1})-\Psi^{n+1}\right)_{K,\sigma}= \\
 \drNnp\left(\frac{Dh(N^{n+1})_{K,\sigma}}{\drNnp}-\frac{D\Psi_{K,\sigma}^{n+1}}{\drNnp}\right).
\end{gather*}
Using the fact that $x=B(-x)-B(x)$ for all $x\in\mathbb{R}$, it rewrites:
\begin{gather*}
D\left(h(N^{n+1})-\Psi^{n+1}\right)_{K,\sigma}=\\ \drNnp\left[B\left(\frac{-Dh(N^{n+1})_{K,\sigma}}{\drNnp}\right)
-B\left(\frac{Dh(N^{n+1})_{K,\sigma}}{\drNnp}\right)\right.\\
 \left. -B\left(\frac{-D\Psi_{K,\sigma}^{n+1}}{\drNnp}\right)
  +B\left(\frac{D\Psi_{K,\sigma}^{n+1}}{\drNnp}\right)\right].
\end{gather*}
But, the definition \eqref{def-dr} of $dr$ implies:
$$ 
\frac{Dh(N^{n+1})_{K,\sigma}}{\drNnp}=D\log(N^{n+1})_{K,\sigma}
$$
and the definition  \eqref{Bern} of $B$ ensures:
\begin{align*}
&B\left(-D\log(N^{n+1})_{K,\sigma}\right)\NKnp-B\left(D\log(N^{n+1})_{K,\sigma}\right)N_{K,\sigma}^{n+1}=0.
\end{align*}
Therefore, we obtain: 
\begin{align*}
\mathcal{R}_{K,\sigma}^{n+1}&=\ts D\left(h(N^{n+1})-\Psi^{n+1}\right)_{K,\sigma}\drNnp\\
&\,\times\left[\left(B\left(\frac{-D\Psi_{K,\sigma}^{n+1}}{\drNnp}\right)-B\left(\frac{-Dh(N^{n+1})_{K,\sigma}}{\drNnp}\right)\right)\right.\\
&\quad\quad\quad\quad\quad\quad\times \left(\NKnp-\min(\NKnp,N_{K,\sigma}^{n+1})\right)\\
&\left.\,\,\,\,\,-\left(B\left(\frac{D\Psi_{K,\sigma}^{n+1}}{\drNnp}\right)-B\left(\frac{Dh(N^{n+1})_{K,\sigma}}{\drNnp}\right)\right)\right. \\
&\quad\quad\quad\quad\quad\quad\times\left.\left(N_{K,\sigma}^{n+1}-\min(\NKnp,N_{K,\sigma}^{n+1})\right)\vphantom{B\left(\frac{-D\Psi_{K,\sigma}^{n+1}}{\drNnp}\right)}\right].
\end{align*}
Since $B$ is nonincreasing on $\mathbb{R}$ and $dr(a,b)\geq 0$ for all $a,\,b\geq 0$, we conclude that $\mathcal{R}_{K,\sigma}^{n+1}\leq 0$. 
\end{proof}


\section{Exponential decay to the discrete thermal equilibrium}\label{sec-exp-decay}

In this section, we establish the main result of this article, namely the decay rate of approximate solutions given by the Scharfetter--Gummel scheme towards an approximation of the thermal equilibrium. Assumptions concerning existence and $L^{\infty}$ estimates will be discussed in the next section. The main theorem is the following:

\begin{theorem}[Exponential decay]\label{thrm:decayE}
Let Hypotheses \ref{HYP} be fulfilled with $m>0$. Let $\M=(\T,\E,\P)$ be an admissible mesh of $\Omega$ satisfying \eqref{regmesh} and $\Delta t$ be the time step verifying \eqref{hyp-dt}. We also assume the compatibility condition \eqref{hyp-compatibility-BC} and that either Hypotheses \ref{HYP-LIN} or Hypotheses \ref{HYP-NON-LIN} are fulfilled.
If moreover there exists a solution $(N_{\T}^{n},P_{\T}^{n},\Psi_{\T}^{n})_{n\geq 0}$ to the numerical scheme \eqref{CIdiscret}--\eqref{scheme-Psi} satisfying
\begin{equation}\label{hyp-borne-unif}
0<m\leq N_{K}^{n},\,P_{K}^{n}\leq M, \quad \forall K\in\T, \quad \forall n\geq 0,
\end{equation}
then there exists a constant $\beta$ only depending on $m$, $M$, $\lambda$, $\Omega$, the function $r$, $\Delta t_{max}$ and the regularity parameter of the mesh $\xi$, but not on the size of the mesh, such that for all $n\geq 0$,
\begin{equation}\label{cv-exp-Ediscret}
\mathbb{E}^{n}\leq e^{-\beta\,t^n}\mathbb{E}^0.
\end{equation}
Furthermore, there exists a constant $c>0$ only depending on $r$, $m$, $\Omega$ and $\xi$ such that for all $n\geq 0$,
\begin{equation}\label{cv-exp-L2-lin}
\|N_{\T}^n-N_{\T}^{eq}\|_{0}^2+\|P_{\T}^n-P_{\T}^{eq}\|_{0}^2+ \|\Psi_{\T}^n-\Psi_{\T}^{eq}\|_{0}^2\leq c\,e^{-\beta\,t^n}\mathbb{E}^0.
 \end{equation} 
\end{theorem}

We notice that the convergence of $\mathbb{E}^{n}$ is exponential, which is in agreement with the exponential decay proved in the continuous framework for the drift--diffusion model both in linear and nonlinear cases \cite{Gajewski1986,Gajewski1989,Gajewski1996}, and for more general electro--reaction--diffusion problems \cite{Gajewski1996a,Glitzky1996,Glitzky2009}. The proof of Theorem \ref{thrm:decayE} is split into two steps. We first establish a discrete analog of \eqref{ineg_EI_continu_1} in Proposition \ref{prop:ineg_EI_discret_1}, and then we prove a control of the relative entropy by the entropy production functional in Proposition \ref{prop:ineg_EI_discret_2}.


\subsection{A first discrete entropy--entropy production inequality}

In the following proposition, we establish the discrete counterpart of \eqref{ineg_EI_continu_1}:

\begin{proposition}\label{prop:ineg_EI_discret_1}
Let Hypotheses \ref{HYP}, the compatibility condition \eqref{hyp-compatibility-BC} and either  Hypotheses \ref{HYP-LIN} or Hypotheses \ref{HYP-NON-LIN} be fulfilled. Let $\M=(\T,\E,\P)$ be an admissible mesh of $\Omega$.
Moreover, if there exists a solution $(N_{\T}^{n},P_{\T}^{n},\Psi_{\T}^{n})_{n\geq 0}$ to \eqref{CIdiscret}--\eqref{scheme-Psi} such that
\begin{equation}\label{hyp-positivite}
0\leq N_{K}^{n},\,P_{K}^{n}\leq M, \quad \forall K\in\T, \quad \forall n\geq 0,
\end{equation}
then for all $n\geq 0$:
\begin{equation}\label{ineg_EI_discret_1}
0\leq \mathbb{E}^{n+1}+\Delta t\,\mathbb{I}^{n+1}\leq \mathbb{E}^{n}.
\end{equation}
\end{proposition}

\begin{proof}
The proof of the isothermal case $r=Id$ with $R=0$ is done in~\cite{Chatard2011}. Here we extend it to the more general case of a nonlinear pressure $r$ satisfying assumptions \eqref{HYP-NON-LIN-1}, or to the linear case with recombination--generation rate of the form \eqref{HYP-LIN-2}.

Using the convexity of $H$ and the definitions \eqref{scheme-N-eq}, \eqref{scheme-P-eq} of $\NKeq$, $\PKeq$, we have:
\begin{multline}\label{estim-prop1-1}
\mathbb{E}^{n+1}-\mathbb{E}^{n}\leq   \somK\mK\left(h(\NKnp)-\alpha_{N}-\Psi_{K}^{eq}\right)\left(\NKnp-\NKn\right)\\
+\somK\mK\left(h(\PKnp)-\alpha_{P}+\Psi_{K}^{eq}\right)\left(\PKnp-\PKn\right)\\
+\frac{\lambda^{2}}{2}\somsige\ts\left[\left(D_{\sigma}\left(\Psi^{n+1}-\Psi^{eq}\right)
\right)^2-\left(D_{\sigma}\left(\Psi^{n}-\Psi^{eq}\right)
\right)^2\right].
\end{multline} 
Since $(a^2-b^2)\leq 2(a-b)a$ for all $a,\,b\in\mathbb{R}$, the third term can be estimated in the following way:
\begin{multline*}
\frac{\lambda^{2}}{2}\somsige\ts\left[\left(D_{\sigma}\left(\Psi^{n+1}-\Psi^{eq}\right)
\right)^2-\left(D_{\sigma}\left(\Psi^{n}-\Psi^{eq}\right)
\right)^2\right]\\
\leq \lambda^2\somsig\!\!\!\ts\,D\left(\Psi^{n+1}-\Psi^{eq}\right)_{K,\sigma}\left[D\left(\Psi^{n+1}-\Psi^{eq}\right)_{K,\sigma}-D\left(\Psi^{n}-\Psi^{eq}\right)_{K,\sigma}\right]\\
\leq \lambda^2\somsig\ts\,D\left(\Psi^{n+1}-\Psi^{eq}\right)_{K,\sigma}\,D\left(\Psi^{n+1}-\Psi^{n}\right)_{K,\sigma}.
\end{multline*}
Then performing a discrete integration by parts and using  \eqref{scheme-Psi} at times $t^n$ and $t^{n+1}$, we get:
\begin{multline*}
\frac{\lambda^{2}}{2}\somsige\ts\left[\left(D_{\sigma}\left(\Psi^{n+1}-\Psi^{eq}\right)
\right)^2-\left(D_{\sigma}\left(\Psi^{n}-\Psi^{eq}\right)
\right)^2\right]\\
\leq -\lambda^2\somK\left(\Psi_{K}^{n+1}-\Psi_{K}^{eq}\right)\sum_{\sigma\in\EK}\ts\,D\left(\Psi^{n+1}-\Psi^{n}\right)_{K,\sigma}\\
\leq \somK\mK\left(\Psi_{K}^{n+1}-\Psi_{K}^{eq}\right)\left[\left(\PKnp-\PKn\right)-\left(\NKnp-\NKn\right)\right].
\end{multline*}
Going back to \eqref{estim-prop1-1}, we obtain:
\begin{multline*}
\mathbb{E}^{n+1}-\mathbb{E}^n\leq  \somK\mK\left(\NKnp-\NKn\right)\left[h(\NKnp)-\Psi_{K}^{n+1}
-\alpha_{N}\right]\\
+\somK\mK\left(\PKnp-\PKn\right)\left[h(\PKnp)+\Psi_{K}^{n+1}-\alpha_{P}
\right].
\end{multline*}
Now using the schemes \eqref{scheme-N} and \eqref{scheme-P}, we get
\begin{equation*}
\mathbb{E}^{n+1}-\mathbb{E}^{n}\leq T_{1}+T_{2}+T_{3},
\end{equation*}
with
\begin{align*}
T_{1}=&-\Delta t\somK\left(h(\NKnp)-\Psi_{K}^{n+1}-\alpha_{N}\right)\sum_{\sigma\in\EK}\mathcal{F}_{K,\sigma}^{n+1},\\
T_{2}=&-\Delta t\somK\left(h(\PKnp)+\Psi_{K}^{n+1}-\alpha_{P}\right)\sum_{\sigma\in\EK}\mathcal{G}_{K,\sigma}^{n+1},\\
T_{3}=&-\Delta t\somK\mK R(\NKnp,\PKnp)\\
&\quad\quad\quad\quad\times\left[h\left(\NKnp\right)+h\left(\PKnp\right)-
h\left(\NKeq\right)-h\left(\PKeq\right)\right].
\end{align*}
The term $T_{3}$ is exactly the last term of the entropy production $\mathbb{I}^n$ defined by \eqref{defIdiscret}. Then integrating $T_{1}$ by parts (since $h(N_{\sigma}^{D})-\Psi^{D}_{\sigma}=\alpha_{N}$ by assumption~\eqref{hyp-compatibility-BC}), we have:
\begin{equation*}
T_{1}=\Delta t\somsig \mathcal{F}_{K,\sigma}^{n+1}D\left(h(N^{n+1})-\Psi^{n+1}\right)_{K,\sigma},
\end{equation*}
and using \eqref{prop-FLUX-N} we obtain
\begin{equation*}
T_{1}\leq -\Delta t\somsig \ts\,\min(\NKnp,N_{K,\sigma}^{n+1})\left(D_{\sigma}\left(h(N^{n+1})-\Psi^{n+1}\right)\right)^{2}.
\end{equation*}
We proceed exactly in the same way for $T_{2}$ by using \eqref{prop-FLUX-P} and get:
$$T_{2}\leq -\Delta t\somsig\ts\,\min(\PKnp,P_{K,\sigma}^{n+1})\left(D_{\sigma}\left(h(P^{n+1})+\Psi^{n+1}\right)\right)^{2},$$
which completes the proof of Proposition \ref{prop:ineg_EI_discret_1}.
\end{proof}


\subsection{Control of the relative entropy by the entropy production}

The proof of Theorem \ref{thrm:decayE} relies on Proposition \ref{prop:ineg_EI_discret_1} and on the following result, which gives a control of the relative entropy by the entropy production:

\begin{proposition}\label{prop:ineg_EI_discret_2}
Let Hypotheses \ref{HYP}, the compatibility condition \eqref{hyp-compatibility-BC} and either  Hypotheses \ref{HYP-LIN} or Hypotheses \ref{HYP-NON-LIN} be fulfilled. We also assume that $m>0$. Let $\M=(\T,\E,\P)$ be an admissible mesh of $\Omega$ satisfying \eqref{regmesh}.
If moreover there exists a solution $(N_{\T}^{n},P_{\T}^{n},\Psi_{\T}^{n})_{n\geq 0}$ to \eqref{CIdiscret}--\eqref{scheme-Psi} such that
\begin{equation}\label{hyp-prop-2}
0< m \leq \NKn,\,\PKn\leq M,\quad \forall K\in\T, \quad \forall n\geq 0,
 \end{equation} 
there exist constants $C_{EF},\,C_{EI}>0$ only depending on $\Omega$, $\xi$, $r$, $\lambda$, $m$ and $M$ such that:
\begin{equation}\label{ineg_EI_discret_2}
\mathbb{E}^{n}\leq C_{EF}\mathbb{F}^n\leq C_{EI}\,\mathbb{I}^{n},\quad  \forall n\geq 0,
\end{equation}
where
\begin{equation*}
\mathbb{F}^n:=\|N^{n}_{\T}-N^{eq}_{\T}\|_{0}^{2}+\|P^{n}_{\T}-P^{eq}_{\T}\|_{0}^{2}+\frac{\lambda^{2}}{2}|\Psi^{n}_{\M}-\Psi_{\M}^{eq}|^{2}_{1,\M},\quad \forall n\geq 0.
\end{equation*}
\end{proposition}

In the proof of Proposition \ref{prop:ineg_EI_discret_2}, we will use the following lemma, whose proof is straightforward  using the assumptions on $r$ and Taylor expansions of $H$ or $h$.
\begin{lemma}
We assume that either $r=Id$ or $r$ satisfies \eqref{HYP-NON-LIN-1}. Let $h$ be the corresponding enthalpy function and $H$ an antiderivative of $h$.
Let $0<m\leq M$. There exist three constants $c_{1}$, $c_{2}$, $c_{3}>0$ depending on $m$, $M$ and $r$ such that for all $x,\,y\in[m,M]$,
\begin{gather}
c_{1}(x-y)^2\leq H(x)-H(y)-h(y)(x-y)\leq c_{2}(x-y)^2, \label{ineg-h-1} \\
c_{3}(x-y)^2\leq(h(x)-h(y))(x-y).\label{ineg-h-2}
\end{gather}
\end{lemma}


Now we can proceed with the proof of Proposition \ref{prop:ineg_EI_discret_2}.

\begin{proof}[Proof of Proposition \ref{prop:ineg_EI_discret_2}]
Let $n\geq 0$. On the one hand, it is clear that 
\begin{equation}\label{estim-E-F}
\mathbb{E}^{n}\leq C_{EF}\,\mathbb{F}^{n},
\end{equation}
using the definition of $\mathbb{E}^n$, estimate \eqref{ineg-h-1} and assumption \eqref{hyp-prop-2}, with $C_{EF}$ only depending on $m$, $M$ and $r$. On the other hand, we prove that 
\begin{equation}\label{estim-F-I}
\mathbb{F}^{n}\leq C_{FI}\,\mathbb{I}^{n}, 
\end{equation}
 by following the same strategy as in the proof of \cite[Theorem 5.3]{Juengel1995}. Using Cauchy--Schwarz and Young inequalities, together with $\alpha_{N}=h(N^{eq})-\Psi^{eq}$, $\alpha_{P}=h(P^{eq})+\Psi^{eq}$, we have for a $\delta>0$, which will be determined later, that:
\begin{align*}
\frac{\delta}{2}&\|N^{n}_{\T}-N_{\T}^{eq}\|_{0}^{2}+\frac{1}{2\,\delta}\|h(N^{n}_{\T})-\Psi^{n}_{\T}-\alpha_{N}\|^{2}_{2}\\
+\frac{\delta}{2}&\|P^{n}_{\T}-P^{eq}_{\T}\|_{0}^{2}+\frac{1}{2\,\delta}\|h(P^{n}_{\T})+\Psi_{\T}^{n}-\alpha_{P}\|^{2}_{2}\\
& \geq \somK\mK\left(\NKn-\NKeq\right)\left(h(\NKn)-\Psi_{K}^{n}-\alpha_{N}\right)
\\&\quad+\somK\mK\left(\PKn-\PKeq
\right)\left(h(\PKn)-\Psi_{K}^{n}-\alpha_{P}\right)\\
&\geq  \somK\mK\left(\NKn-\NKeq\right)\left(h(\NKn)-h(\NKeq)-\left(\Psi_{K}^{n}-\Psi_{K}^{eq}\right)\right)\\
& \quad\quad +\somK\mK\left(\PKn-\PKeq\right)\left(h(\PKn)-h(\PKeq)+\left(\Psi_{K}^{n}-\Psi_{K}^{eq}\right)\right).
\end{align*} 
Using estimate \eqref{ineg-h-2} and the schemes \eqref{scheme-Psi} and \eqref{scheme-Psi-eq}, we obtain by integrating by parts ($\Psi_{\sigma}^{n}=\Psi_{\sigma}^{D}=\Psi_{\sigma}^{eq}$ for all $\sigma\in\E_{ext}^{D}$) that
\begin{align*}
\frac{\delta}{2}&\|N^{n}_{\T}-N_{\T}^{eq}\|_{0}^{2}+\frac{1}{2\,\delta}\|h(N_{\T}^{n})-\Psi_{\T}^{n}-\alpha_{N}\|^{2}_{2}\\
+\frac{\delta}{2}&\|P_{\T}^{n}-P_{\T}^{eq}\|_{0}^{2}+\frac{1}{2\,\delta}\|h(P_{\T}^{n})+\Psi_{\T}^{n}-\alpha_{P}\|^{2}_{2}\\
& \geq c_{3}\|N^{n}_{\T}-N^{eq}_{\T}\|_{0}^{2}+c_{3}\|P^{n}_{\T}-P^{eq}_{\T}\|_{0}^{2}\\
&\quad\quad+\somK\mK\left[\left(\PKn-\NKn+C_{K}\right)-\left(\PKeq-\NKeq+C_{K}\right)\right](\Psi_{K}^{n}-\Psi_{K}^{eq})\\
&\geq  c_{3}\|N^{n}_{\T}-N^{eq}_{\T}\|_{0}^{2}+c_{3}\|P^{n}_{\T}-P^{eq}_{\T}\|_{0}^{2}\\
&\quad\quad-\lambda^{2}\somK\sum_{\sigma\in\EK}\ts\,D\left(\Psi^{n}-\Psi^{eq}\right)_{K,\sigma}(\Psi_{K}^{n}-\Psi_{K}^{eq})\\
&\geq  c_{3}\|N^{n}_{\T}-N^{eq}_{\T}\|_{0}^{2}+c_{3}\|P^{n}_{\T}-P^{eq}_{\T}\|_{0}^{2}+\lambda^{2}|\Psi_{\M}^{n}-\Psi_{\M}^{eq}|_{1,\M}^{2}.
\end{align*}
Then we obtain that
\begin{gather*}
\left(c_{3}-\frac{\delta}{2}\right)\left[\|N^{n}_{\T}-N^{eq}_{\T}\|_{0}^{2}+\|P^{n}_{\T}-P^{eq}_{\T}\|_{0}^{2}\right]+\lambda^{2}|\Psi_{\M}^{n}-\Psi_{\M}^{eq}|^{2}_{1,\M}\\
\leq \frac{1}{2\,\delta}\|h(N^{n}_{\T})-\Psi_{\T}^{n}-\alpha_{N}\|_{0}^{2}+\frac{1}{2\,\delta}\|h(P^{n}_{\T})+\Psi_{\T}^{n}-\alpha_{P}\|_{0}^{2}.
\end{gather*}
Now choosing $\delta=c_{3}$ and applying the discrete Poincar\'e inequality  \eqref{ineq_Poincare} (since $h(N^{n}_{\sigma})-\Psi_{\sigma}^{n}-\alpha_{N}=0$ and $h(P^{n}_{\sigma})+\Psi_{\sigma}^{n}-\alpha_{P}=0$ for all $\sigma\in\E_{ext}^{D}$ and $\text{m}(\Gamma^{D})>0$), it yields the existence of a constant $c>0$ depending only on $m$, $r$, $\lambda$, $\xi$, $\Omega$ such that:
\begin{equation*}
\mathbb{F}^{n}\leq c\left(|h(N^{n}_{\M})-\Psi^{n}_{\M}|^{2}_{1,\M}+|h(P^{n}_{\M})+\Psi^{n}_{\M}|^{2}_{1,\M}\right).
\end{equation*}
But by definition $|\cdot|_{1,\M}$ and using the uniform lower bound of $N_{\M}$, $P_{\M}$, we have:
\begin{gather*}
 |h(N^{n}_{\M})-\Psi^{n}_{\M}|^{2}_{1,\M}+|h(P^{n}_{\M})+\Psi^{n}_{\M}|^{2}_{1,\M}\leq\\
 \frac{1}{m}\somsig\ts\,\min(\NKn,N_{K,\sigma}^{n})\left(D_{\sigma}\left(h(N^{n})-\Psi^{n}\right)\right)^{2}\\
 +\frac{1}{m}\somsig\ts\,\min(\PKn,P_{K,\sigma}^{n})\left(D_{\sigma}\left(h(P^{n})+\Psi^{n}\right)\right)^{2}\\
 \leq\frac{1}{m}\,\mathbb{I}^{n},
\end{gather*}
which concludes the proof of \eqref{estim-F-I}.
\end{proof}


\subsection{Proof of Theorem \ref{thrm:decayE}}

In view of Propositions \ref{prop:ineg_EI_discret_1} and \ref{prop:ineg_EI_discret_2}, it is now easy to prove Theorem \ref{thrm:decayE}. Indeed, we have
\begin{equation*}
\mathbb{E}^{n+1}-\mathbb{E}^{n}\leq -\Delta t\,\mathbb{I}^{n+1}\leq -\frac{\Delta t}{C_{EI}}\mathbb{E}^{n+1}.
\end{equation*}
We use this inequality to prove that there exists $\beta>0 $ such that $\mathbb{E}^n\leq e^{-\beta t^n}\mathbb{E}^0$. Indeed, we have
\begin{align*}
e^{\beta t^{n+1}}\mathbb{E}^{n+1}-e^{\beta t^{n}}\mathbb{E}^{n} & =\left(e^{\beta t^{n+1}}-e^{\beta t^{n}}\right)\mathbb{E}^{n+1}+e^{\beta t^{n}}\left(\mathbb{E}^{n+1}-\mathbb{E}^n\right)\\
&\leq e^{\beta t^n}\Delta t\left[e^{\theta \Delta t}\beta-\frac{1}{C_{EI}}\right]\mathbb{E}^{n+1},
\end{align*}
where $\theta\in (0,1)$. Then using assumption \eqref{hyp-dt} and choosing $\beta=e^{-\Delta t_{max}}/C_{EI}$, we have
$$e^{\theta\Delta t}\beta-\frac{1}{C_{EI}}\leq 0.$$
Thus the sequence $\left(e^{\beta t^n}\mathbb{E}^n\right)_{n}$ is nonincreasing, which yields the result.\\
Then the $L^2$ estimate \eqref{cv-exp-L2-lin} is straightforward using that inequality \eqref{ineg-h-1} and the discrete Poincar\'e inequality \eqref{ineq_Poincare} imply
$$\|N_{\T}^n-N_{\T}^{eq}\|_{0}^2+\|P_{\T}^n-P_{\T}^{eq}\|_{0}^2+ \|\Psi_{\T}^n-\Psi_{\T}^{eq}\|_{0}^2 \leq C \mathbb{E}^n.$$


\section{Existence of a numerical solution and $L^{\infty}$ estimates}\label{sec-existence}


In this section, we discuss about the assumptions made in Theorem \ref{thrm:decayE} concerning existence of a solution to the numerical scheme and uniform $L^{\infty}$ estimates. We distinguish the isothermal case satisfying Hypotheses \ref{HYP-LIN} from the nonlinear case without recombination--generation rate satisfying Hypotheses~\ref{HYP-NON-LIN}.


\subsection{The nonlinear case without recombination--generation rate}

\begin{theorem}\label{thrm-existence-nonlin}
Let Hypotheses \ref{HYP} and Hypotheses \ref{HYP-NON-LIN} be fulfilled. Let $\M=(\T,\E,\P)$ be an admissible mesh of $\Omega$. Moreover we assume that the time step satisfies:
\begin{equation}\label{hyp-dt-nonlin}
\Delta t\leq \frac{\lambda^2}{\|C\|_{\infty}}.
\end{equation}
Then there exists a solution $(N_{\T}^{n},P_{\T}^{n},\Psi_{\T}^{n})_{n\geq 0}$ to \eqref{CIdiscret}--\eqref{scheme-Psi}, which  satisfies the following $L^{\infty}$ estimates for the approximate densities: $\forall K\in\T, \quad \forall 0\leq n\leq T/\Delta t,$
\begin{equation}\label{borne-Linf-gene}
m\exp\left(-\frac{\|C\|_{\infty}}{\lambda^2}T\right)\leq m^n \leq \NKn,\,\PKn\leq M^n\leq M\exp\left(\frac{\|C\|_{\infty}}{\lambda^2}T\right),
\end{equation} 
where
\begin{equation*}
m^n=m\left(1+\frac{\Delta t}{\lambda^2}\|C\|_{\infty}\right)^{-n},\quad M^n=M\left(1-\frac{\Delta t}{\lambda^2}\|C\|_{\infty}\right)^{-n}.
\end{equation*}
In particular, if $C=0$, the maximum principle holds for the densities:
\begin{equation}\label{borne-unif}
\forall K\in\T,\quad\forall n\geq 0,\quad m\leq \NKn,\,\PKn\leq M.
\end{equation}
\end{theorem}

Let us emphasize that in the zero doping case, we get the existence of a solution satisfying the assumption \eqref{hyp-borne-unif} of Theorem \ref{thrm:decayE} (providing $m>0$), without any restricting assumption on the time step.

\begin{proof}
We prove the result by induction on $n\geq 0$. The vectors $N_{\T}^{0}$ and $P_{\T}^{0}$ are given by \eqref{CIdiscret} while $\Psi_{\T}^{0}$ is uniquely defined by \eqref{scheme-Psi}. Then, the assumption \eqref{HYP-M-CI-BC} on the initial data  ensures that 
$$m \leq N_{K}^{0},\,P_{K}^{0}\leq M\quad\forall K\in\T.
$$
Now we suppose that, for some $n\geq 0$, $(N_{\T}^{n},P_{\T}^{n},\Psi_{\T}^{n})$ is known and satisfies the $L^{\infty}$ estimate~\eqref{borne-Linf-gene}. We have to establish the existence of $(N_{\T}^{n+1},P_{\T}^{n+1},\Psi_{\T}^{n+1})$ solution to the nonlinear system of equations \eqref{scheme-N}--\eqref{scheme-Psi} satisfying \eqref{borne-Linf-gene} with $n+1$ instead of $n$. We extend the proof done in \cite{Bessemoulin-Chatard2014a}, which follows some ideas developed in \cite{Prohl2009}, to the nonlinear case with nonvanishing doping profile. The method consists in introducing a problem penalized by an arbitrary parameter which will be conveniently chosen.\\
Let $\mu>0$. We introduce an application $T_{\mu}^{n}:\mathbb{R}^{\theta}\times\mathbb{R}^{\theta}\rightarrow \mathbb{R}^{\theta}\times\mathbb{R}^{\theta}$ such that $T_{\mu}^{n}(N_{\T},P_{\T})=(\hat{N}_{\T},\hat{P}_{\T})$, based on a linearization of the scheme \eqref{scheme-N}--\eqref{scheme-Psi} and defined in two steps.
\begin{itemize}
\item Step 1: we define $\Psi_{\T}\in\mathbb{R}^{\theta}$ as the solution to the following linear system:
\begin{equation}\label{scheme-Psi-lin}
-\lambda^{2}\sum_{\sigma\in\E_{K}}\ts\,D\Psi_{K,\sigma}=\mK(P_{K}-N_{K}+C_{K})\quad\forall K\in \T,
\end{equation}
with $\Psi_{\sigma}=\Psi_{\sigma}^{D}$ for all $\sigma\in\E_{ext}^{D}$.
\item Step 2: we construct $(\hat{N}_{\T},\hat{P}_{\T})$ as the solution to the following linear scheme: for all $K\in\T$,
\begin{gather}
\frac{\mK}{\Delta t}\left[\left(1+\frac{\mu}{\lambda^{2}}\right)\hat{N}_{K}-\frac{\mu}{\lambda^{2}}N_{K}-\NKn\right]
+\sum_{\sigma\in\E_{K}}\ts\,dr(N_{K},N_{K,\sigma})\label{scheme-N-lin}\\
\times\left[B\left(\frac{-D\Psi_{K,\sigma}}{dr(N_{K},N_{K,\sigma})}\right)\hat{N}_{K}-B\left(\frac{D\Psi_{K,\sigma}}{dr(N_{K},N_{K,\sigma})}\right)\hat{N}_{K,\sigma}\right]
=0,\nonumber
\end{gather}
\begin{gather}
\frac{\mK}{\Delta t}\left[\left(1+\frac{\mu}{\lambda^{2}}\right)\hat{P}_{K}-\frac{\mu}{\lambda^{2}}P_{K}-\PKn\right]
+\sum_{\sigma\in\E_{K}}\ts\,dr(P_{K},P_{K,\sigma})\label{scheme-P-lin}\\
\times\left[B\left(\frac{D\Psi_{K,\sigma}}{dr(P_{K},P_{K,\sigma})}\right)\hat{P}_{K}-B\left(\frac{-D\Psi_{K,\sigma}}{dr(P_{K},P_{K,\sigma})}\right)\hat{P}_{K,\sigma}\right]=0,\nonumber
\end{gather}
with $\hat{N}_{\sigma}=N_{\sigma}^{D}$ and $\hat{P}_{\sigma}=P_{\sigma}^{D}$ for all $\sigma\in\E_{ext}^{D}$.
\end{itemize}
The existence and uniqueness of $\Psi_{\T}$ solution to the linear system \eqref{scheme-Psi-lin} are obvious. Schemes \eqref{scheme-N-lin} and \eqref{scheme-P-lin} also lead to two decoupled linear systems which can be written under a matricial form: $\mathbb{A}_{N}\hat{N}_{\T}=\mathbb{S}_{N}^{n}$ and $\mathbb{A}_{P}\hat{P}_{\T}=\mathbb{S}_{P}^{n}$. The matrix $\mathbb{A}_{N}$ is the sparse matrix defined by: $\forall K \in\T,$
\begin{align*}
&(\mathbb{A}_{N})_{K,K}=\frac{\mK}{\Delta t}\left(1+\frac{\mu}{\lambda^{2}}\right)+\!\!\!\!\!\!\sum_{\sigma\in\EK\setminus\E_{K,ext}^{N}}\!\!\!\!\ts\,dr(N_{K},N_{K,\sigma})B\left(\frac{-D\Psi_{K,\sigma}}{dr(N_{K},N_{K,\sigma})}\right),\\
&(\mathbb{A}_{N})_{K,L}=-\ts\,dr(N_{K},N_{L})B\left(\frac{D\Psi_{K,\sigma}}{dr(N_{K},N_{L})}\right),\,  \forall L\in\T \text{ s. t. } \sigma=K|L\in\Eint.
\end{align*}
The matrix $\mathbb{A}_{N}$ has positive diagonal terms, nonpositive offdiagonal terms and is strictly diagonally dominant with respect to its columns. Then $\mathbb{A}_{N}$ is an M-matrix, which implies that it is invertible and its inverse has only nonnegative coefficients. The same result holds for $\mathbb{A}_{P}$. Thus we obtain that the scheme \eqref{scheme-N-lin}--\eqref{scheme-P-lin} admits a unique solution $(\hat{N}_{\T},\hat{P}_{\T})\in\mathbb{R}^{\theta}\times\mathbb{R}^{\theta}$, so that the application $T_{\mu}^{n}$ is well-defined and is moreover continuous.\\
Now in order to apply the Brouwer's fixed point theorem, we prove that $T_{\mu}^{n}$ preserves the set
\begin{equation*}
\mathcal{C}_{n+1}=\{(N_{\T},P_{\T})\in\mathbb{R}^{\theta}\times\mathbb{R}^{\theta};\quad m^{n+1}\leq N_{K},\,P_{K}\leq M^{n+1}\}.
\end{equation*}
The right hand side of the system $\mathbb{A}_{N}\hat{N}_{\T}=\mathbb{S}_{N}^{n}$ is defined by
\begin{gather*}
(\mathbb{S}_{N}^{n})_{K}=\frac{\mK}{\Delta t}\left(\frac{\mu}{\lambda^{2}}N_{K}+\NKn\right)\\
+\sum_{\sigma\in\E_{K,ext}^{D}}\ts\,dr(N_{K},N_{\sigma}^{D})B\left(\frac{
D\Psi_{K,\sigma}}{dr(N_{K},N_{\sigma}^{D})}\right)N_{\sigma}^{D} \quad \forall K\in\T.
\end{gather*}
If $N_{\T}\geq 0$, then $\mathbb{S}^{n}_{N}\geq 0$ and since $\mathbb{A}_{N}$ is an M-matrix, we obtain that $\hat{N}_{\T}\geq 0$. In the same way, if $P_{\T}\geq 0$, we obtain that $\hat{P}_{\T}\geq 0$.

In order to prove that $\hat{N}_{K}\leq M^{n+1}$ for all $K\in\T$, we introduce now $\mathbf{M}^{n+1}_{\T}$, the constant vector of $\mathbb{R}^{\theta}$ with unique value $M^{n+1}$, and we compute $\mathbb{A}_{N}(\hat{N}_{\T}-\mathbf{M}^{n+1}_{\T})$. For all $K\in\T$, using the fact that $B(x)-B(-x)=-x$ for all $x\in\mathbb{R}$, we have
\begin{multline*}
\left(\mathbb{A}_{N}(\hat{N}_{\T}-\mathbf{M}^{n+1}_{\T})\right)_{K}=\frac{\mK}{\Delta t}\left(\NKn-M^{n+1}\right)+\frac{\mK}{\Delta t}\frac{\mu}{\lambda^2}\left(N_{K}-M^{n+1}\right)\\
-\somsigint\ts D\Psi_{K,\sigma}\,M^{n+1}
 +\sum_{\sigma\in\E_{K,ext}^D}\ts dr(N_{K},N_{\sigma}^D)\left(B\left(\frac{D\Psi_{K,\sigma}}{dr(N_{K},N_{\sigma}^D)}\right)N_{\sigma}^D\right.\\
 \left.-B\left(\frac{-D\Psi_{K,\sigma}}{dr(N_{K},N_{\sigma}^D)}\right)M^{n+1}\right).
\end{multline*}
Since $B$ is a nonnegative function and $N_{\sigma}^D\leq M\leq M^{n+1}$ for all $\sigma\in\E_{ext}^D$, we obtain
\begin{multline*}
\left(\mathbb{A}_{N}(\hat{N}_{\T}-\mathbf{M}^{n+1}_{\T})\right)_{K}\leq  \frac{\mK}{\Delta t}\left(\NKn-M^n\right)+\frac{\mK}{\Delta t}\left(M^n-M^{n+1}\right)\\
+\frac{\mK}{\Delta t}\frac{\mu}{\lambda^2}\left(N_{K}-M^{n+1}\right)-\sum_{\sigma\in\E_{K}}\ts D\Psi_{K,\sigma}M^{n+1}.
\end{multline*}
The induction assumption ensures that $\NKn\leq M^n$. Since 
\begin{equation*}
M^n-M^{n+1}=-\ds{\frac{\Delta t}{\lambda^2}\|C\|_{\infty}M^{n+1}},
\end{equation*}
 the scheme \eqref{scheme-Psi-lin} leads to:
\begin{align}
\left(\mathbb{A}_{N}(\hat{N}_{\T}-\mathbf{M}^{n+1}_{\T})\right)_{K}&\leq  -\frac{\mK}{\lambda^2}\|C\|_{\infty}M^{n+1}+\frac{\mK}{\Delta t}\frac{\mu}{\lambda^2}\left(N_{K}-M^{n+1}\right)\nonumber\\
&+\frac{\mK}{\lambda^2}\left((P_{K}-M^{n+1})-(N_{K}-M^{n+1})+C_{K}\right)M^{n+1}\nonumber\\&
\leq  \frac{\mK}{\lambda^2}\left(\frac{\mu}{\Delta t}-M^{n+1}\right)(N_{K}-M^{n+1})\nonumber \\
&+\frac{\mK}{\lambda^2}(P_{K}-M^{n+1})M^{n+1}\nonumber\\
&+\frac{\mK}{\lambda^2}M^{n+1}\left(C_{K}-\|C\|_{\infty}\right).\label{ineg-N-M}
\end{align}
We can prove exactly in the same way that for all $K\in\T$,
\begin{multline}\label{ineg-N-m}
\left(\mathbb{A}_{N}(\hat{N}_{\T}-\mathbf{m}^{n+1}_{\T})\right)_{K}\geq \frac{\mK}{\lambda^2}\left(\frac{\mu}{\Delta t}-m^{n+1}\right)(N_{K}-m^{n+1})\\
+\frac{\mK}{\lambda^2}m^{n+1}(P_{K}-m^{n+1})\\
+\frac{\mK}{\lambda^2}m^{n+1}\left(\|C\|_{\infty}+C_{K}\right).
\end{multline}
Since $\mu>0$ is an arbitrary constant, we can choose it in such a way that $\mu\geq M^{n+1}\Delta t \geq m^{n+1}\Delta t$. Then if $(N_{\T},P_{\T})\in\mathcal{C}_{n+1}$, inequalities \eqref{ineg-N-M} and \eqref{ineg-N-m} imply that 
$$ \mathbb{A}_{N}(\hat{N}_{\T}-\mathbf{M}^{n+1}_{\T})\leq 0\text{ and } \mathbb{A}_{N}(\hat{N}_{\T}-\mathbf{m}^{n+1}_{\T})\geq 0.$$
Since $\mathbb{A}_{N}$ is an M-matrix, we conclude that $m^{n+1}\leq \hat{N}_{K}\leq M^{n+1}$ for all $K\in\T$. The proof that $m^{n+1}\leq\hat{P}_{K}\leq M^{n+1}$ for all $K\in\T$ is similar and then we have $(\hat{N}_{\T},\hat{P}_{\T})\in\mathcal{C}_{n+1}$.

Finally, $T_{\mu}^{n}$ is a continuous application  which stabilizes the set $\mathcal{C}_{n+1}$. Then, by the Brouwer's fixed point theorem, $T_{\mu}^{n}$ admits a fixed point in $\mathcal{C}_{n+1}$, which is denoted by $(N_{\T}^{n+1},P_{\T}^{n+1})$ and satisfies the $L^{\infty}$ estimate \eqref{borne-Linf-gene}. The corresponding $\Psi_{\T}$ defined by \eqref{scheme-Psi-lin} is denoted by $\Psi_{\T}^{n+1}$ and $(N_{\T}^{n+1},P_{\T}^{n+1},\Psi_{\T}^{n+1})$ is a solution to the scheme \eqref{scheme-N}--\eqref{scheme-Psi}, which concludes the proof.
\end{proof}


\subsection{The isothermal case with recombination--generation rate}

\begin{proposition}\label{prop-existence-R}
Let Hypotheses \ref{HYP} and Hypotheses \ref{HYP-LIN} be fulfilled, with $m=1/M>0$ in \eqref{HYP-M-CI-BC}, namely
\begin{equation}\label{hyp-borne-CI-BC-R}
\frac{1}{M}\leq N^D,\,P^D,\,N_{0},\,P_{0}\leq M.
\end{equation}
Let $\M=(\T,\E,\P)$ be an admissible mesh of $\Omega$. We also assume that the time step satisfies \eqref{hyp-dt-nonlin}.
 Then there exists a solution $(N_{\T}^n,P_{\T}^n,\Psi_{\T}^n)_{n\geq 0}$ to the scheme \eqref{CIdiscret}--\eqref{FLUX-SG-P}, which satisfies the following $L^{\infty}$ estimates for the approximate densities: 
\begin{equation}\label{borne-Linf-gene-R}
0\leq \frac{1}{M^n}\leq N_{K}^n,\,P_{K}^n\leq M^n, \quad \forall K\in\T,\ \forall 0\leq n\leq T/\Delta t.
\end{equation}
where $\ds{M^n=M\left(1-\frac{\Delta t}{\lambda^2}\|C\|_{\infty}\right)^{-n} \leq  M\exp\left(\frac{\|C\|_{\infty}}{\lambda^2}T\right) }.$
\end{proposition}

\begin{proof}
We proceed as in the proof of Theorem \ref{thrm-existence-nonlin}, by induction on $n\geq 0$. We suppose that for some $n\geq 0$, $(N_{\T}^n,P_{\T}^n,\Psi_{\T}^n)$ is known and satisfies \eqref{borne-Linf-gene-R}.

Let $\mu>0$. We define $T_{\mu}^n:(N_{\T},P_{\T})\in\mathbb{R}^{\theta}\times \mathbb{R}^{\theta}\mapsto (\hat{N}_{\T},\hat{P}_{\T})\in \mathbb{R}^{\theta}\times \mathbb{R}^{\theta}$ in two steps. We first define $\Psi_{\T}\in\mathbb{R}^{\theta}$ as the solution to \eqref{scheme-Psi-lin}. Then we construct $(\hat{N}_{\T},\hat{P}_{\T})$ as the solution to the following linear scheme: for all $K\in\T$,
\begin{gather}
\frac{\mK}{\Delta t}\left[\left(1+\frac{\mu}{\lambda^2}\right)\hat{N}_{K}-\frac{\mu}{\lambda^2}N_{K}-\NKn\right]\nonumber\\
+\sum_{\sigma\in\E_{K}}\ts\left[B\left(-D\Psi_{K,\sigma}\right)\hat{N}_{K}-B\left(D\Psi_{K,\sigma}\right)\hat{N}_{K,\sigma}\right] \label{scheme-N-lin-R}\\
=-\mK R_{0}(N_{K},P_{K})\left(\hat{N}_{K}P_{K}-1\right),\nonumber
\end{gather}
\begin{gather}
\frac{\mK}{\Delta t}\left[\left(1+\frac{\mu}{\lambda^2}\right)\hat{P}_{K}-\frac{\mu}{\lambda^2}P_{K}-\PKn\right]\nonumber \\
+\sum_{\sigma\in\E_{K}}\ts\left[B\left(D\Psi_{K,\sigma}\right)\hat{P}_{K}-B\left(-D\Psi_{K,\sigma}\right)\hat{P}_{K,\sigma}\right] \label{scheme-P-lin-R}\\
=-\mK R_{0}(N_{K},P_{K})\left(N_{K}\hat{P}_{K}-1\right),\nonumber
\end{gather}
with $\hat{N}_{\sigma}=N_{\sigma}^D$ and $\hat{P}_{\sigma}=P_{\sigma}^D$ for all $\sigma\in\E_{ext}^D$.
Schemes \eqref{scheme-N-lin-R} and \eqref{scheme-P-lin-R} can be written under a matricial form: $\mathbb{A}_{N}\hat{N}_{\T}=\mathbb{S}_{N}^n$ and $\mathbb{A}_{P}\hat{P}_{\T}=\mathbb{S}_{P}^n$, where $\mathbb{A}_{N}$ is the sparse matrix defined by: $\forall K\in\T$,
\begin{align*}
&\left(\mathbb{A}_{N}\right)_{K,K}=\frac{\mK}{\Delta t}\left(1+\frac{\mu}{\lambda^2}\right)\\
&\quad\quad\quad+\sum_{\sigma\in\E_{K}\setminus\E_{K,ext}^N}\ts\, B\left(-D\Psi_{K,\sigma}\right)+\mK R_{0}(N_{K},P_{K})P_{K},\\
&\left(\mathbb{A}_{N}\right)_{K,L}=-\ts \,B\left(D\Psi_{K,\sigma}\right) \quad\forall L\in\T \text{ such that }\sigma=K|L\in\Eint.
\end{align*}
Since $\mathbb{A}_{N}$ and $\mathbb{A}_{P}$ are M-matrices and $R_{0}$ is continuous, the application $T_{\mu}^n$ is well-defined and continuous. We now prove that it preserves the set 
$$\mathcal{C}_{n+1}=\left\{(N_{\T},P_{\T})\in\mathbb{R}^{\theta}\times
\mathbb{R}^{\theta}; \quad \frac{1}{M^{n+1}}\leq N_{K},\,P_{K}\leq M^{n+1}\right\}.$$
The right hand side $\mathbb{S}_{N}^n$ is defined by
\begin{gather*}
\left(\mathbb{S}_{N}^n\right)_{K}=\frac{\mK}{\Delta t}\left(\frac{\mu}{\lambda^2}N_{K}+\NKn\right)\\
+\sum_{\sigma\in\E_{K,ext}^D}\ts\,B\left(D\Psi_{K,\sigma}\right)N_{\sigma}^D+\mK\,R_{0}(N_{K},P_{K}).
\end{gather*}
It is clear that if $N_{\T}\geq 0$, $P_{\T}\geq 0$, then $\hat{N}_{\T}\geq 0$, $\hat{P}_{\T}\geq 0$. Now to prove that $\hat{N}_{K}\leq M^{n+1}$, we compute $\mathbb{A}_{N}\left(\hat{N}_{\T}-\mathbf{M}^{n+1}_{\T}\right)$ as in Section \ref{sec-existence}: for all $K\in\T$, we obtain
\begin{align*}
\left(\mathbb{A}_{N}\left(\hat{N}_{\T}-\mathbf{M}^{n+1}_{\T}\right)
\right)_{K}\leq & \frac{\mK}{\lambda^2}\left(\frac{\mu}{\Delta t}-M^{n+1}\right)\left(N_{K}-M^{n+1}\right)\\
&+\frac{\mK}{\lambda^2}M^{n+1}\left(P_{K}-M^{n+1}\right)\\
&+\frac{\mK}{\lambda^2}M^{n+1}\left(C_{K}-\|C\|_{\infty}\right)\\
&+\mK R_{0}(N_{K},P_{K})\left(1-P_{K}M^{n+1}\right)\\
\leq & \frac{\mK}{\lambda^2}\left(\frac{\mu}{\Delta t}-M^{n+1}\right)\left(N_{K}-M^{n+1}\right).
\end{align*}
Since $(N_{\T},P_{\T})\in\mathcal{C}_{n+1}$ and $R_{0}(N_{K},P_{K})\geq 0$, this yields
\begin{equation}\label{ineg-N-M-R}
\left(\mathbb{A}_{N}\left(\hat{N}_{\T}-\mathbf{M}^{n+1}_{\T}\right)
\right)_{K}\leq \frac{\mK}{\lambda^2}\left(\frac{\mu}{\Delta t}-M^{n+1}\right)\left(N_{K}-M^{n+1}\right).
\end{equation}
If we define $\mathbf{m}^{n+1}$ the constant vector of $\mathbb{R}^{\theta}$ with unique value $\frac{1}{M^{n+1}}$, we obtain
\begin{align*}
\left(\mathbb{A}_{N}\left(\hat{N}_{\T}-\mathbf{m}^{n+1}_{\T}\right)
\right)_{K}\geq & \frac{\mK}{\lambda^2}\left(\frac{\mu}{\Delta t}-\frac{1}{M^{n+1}}\right)\left(N_{K}-\frac{1}{M^{n+1}}\right)\\
& +\frac{\mK}{\lambda^2}\frac{1}{M^{n+1}}\left(P_{K}-\frac{1}{M^{n+1}}\right) \\
&+\mK R_{0}(N_{K},P_{K})\left(1-\frac{P_{K}}{M^{n+1}}\right)\\
&+\frac{\mK}{\Delta t}\left(\frac{1}{M^{n+1}}-\frac{1}{M^n}\right)+\frac{\mK}{\lambda^2}\frac{1}{M^{n+1}}C_{K}.
\end{align*}
But using that $M^{n+1}-M^n=\frac{\Delta t}{\lambda^2}\|C\|_{\infty}M^{n+1}$, we have
\begin{align*}
\frac{\mK}{\Delta t}&\left(\frac{1}{M^{n+1}}-\frac{1}{M^n}\right)+\frac{\mK}{\lambda^2}\frac{1}{M^{n+1}}C_{K} \\
&=\frac{\mK}{\lambda^2}\frac{1}{M^{n+1}}\left(\frac{\lambda^2}{\Delta t}\left(\frac{M^{n+1}-M^n}{M^n}\right)+C_{K}\right)\\
&=\frac{\mK}{\lambda^2}\frac{1}{M^{n+1}}\left(\frac{M^{n+1}}{M^n}\|C\|_{\infty}+C_{K}\right)\\
&\geq \frac{\mK}{\lambda^2}\frac{1}{M^{n+1}}\left(\|C\|_{\infty}+C_{K}\right)\geq 0.
\end{align*}
Finally, since $(N_{\T},P_{\T})\in\mathcal{C}_{n+1}$, we obtain
\begin{equation}\label{ineg-N-m-R}
\left(\mathbb{A}_{N}\left(\hat{N}_{\T}-\mathbf{m}^{n+1}_{\T}\right)
\right)_{K}\geq \frac{\mK}{\lambda^2}\left(\frac{\mu}{\Delta t}-\frac{1}{M^{n+1}}\right)\left(N_{K}-\frac{1}{M^{n+1}}\right).
\end{equation}
Then if we choose $\mu>0$ such that $\mu\geq M^{n+1}\Delta t\geq \Delta t/M^{n+1}$, inequalities \eqref{ineg-N-M-R} and \eqref{ineg-N-m-R} imply that $(\hat{N}_{\T},\hat{P}_{\T})\in\mathcal{C}_{n+1}$, and we conclude as in the proof of Theorem \ref{thrm-existence-nonlin}.
\end{proof}


\section{Numerical experiments}\label{sec-numerical-exp}

\begin{figure}[ht!]
\includegraphics[scale=0.4]{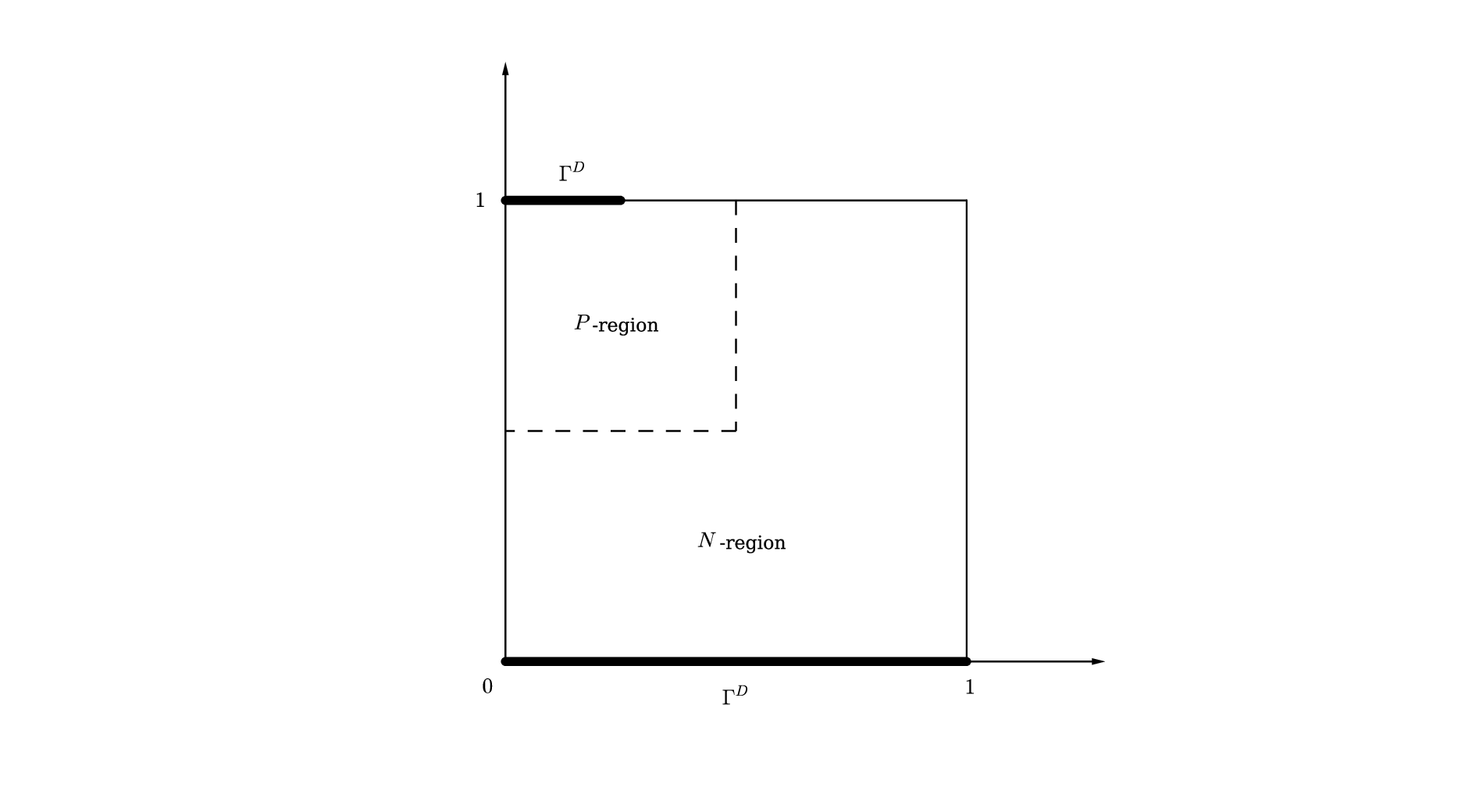}
\caption{Geometry of the PN-junction diode}
\label{jonctionPN}
\end{figure}

In this section, we present some numerical experiments that illustrate the long time behavior of approximate solutions in various situations. We consider a geometry corresponding to a PN-junction in 2D (see Figure \ref{jonctionPN}). The domain $\Omega$ is the square $(0,1)^2$. The Dirichlet boundary conditions are:
\begin{align*}
& N^D=N^D_{0},\,P^D=P^D_{0} \quad \text{ on } \{y=0\},\\
& N^D=N^D_{1},\,P^D=P^D_{1} \quad \text{ on } \{y=1,\,0\leq x\leq 0.25\},
\end{align*}
for values of $N^D_{0}$,$N^D_{1}$, $P^D_{0}$, $P^D_{1}$ such that  $h(N^D)+h(P^D)=c_{0}$ on $\Gamma^D$. Then we define $\Psi^D=\left(h(N^D)-h(P^D)\right)/2$, in such a way that \eqref{hyp-compatibility-BC} is satisfied with $\alpha_{N}=\alpha_{P}=c_{0}/2$. Elsewhere we put homogeneous Neumann boundary conditions. Initial conditions are given by $N_{0}(x,y)=N_{1}^D+(N_{0}^D-N_{1}^D)(1-\sqrt{y})$, $P_{0}(x,y)=P_{1}^D+(P_{0}^D-P_{1}^D)(1-\sqrt{y})$.

We compute the numerical approximation of the thermal equilibrium and of the transient drift--diffusion system on a mesh made of 3584 triangles, with time step $\Delta t=10^{-2}$, until time $T=10$ or $20$. Since we study the convergence towards the equilibrium and not the quasi--neutral limit, we fix $\lambda^2=1$. Finally, the doping profile is either zero or piecewise constant, equal to 1 in the N-region and -1 in the P-region.

\begin{figure}
\centering
\subfigure[$C=0$, $R=0$]{\includegraphics[width=1.6in]{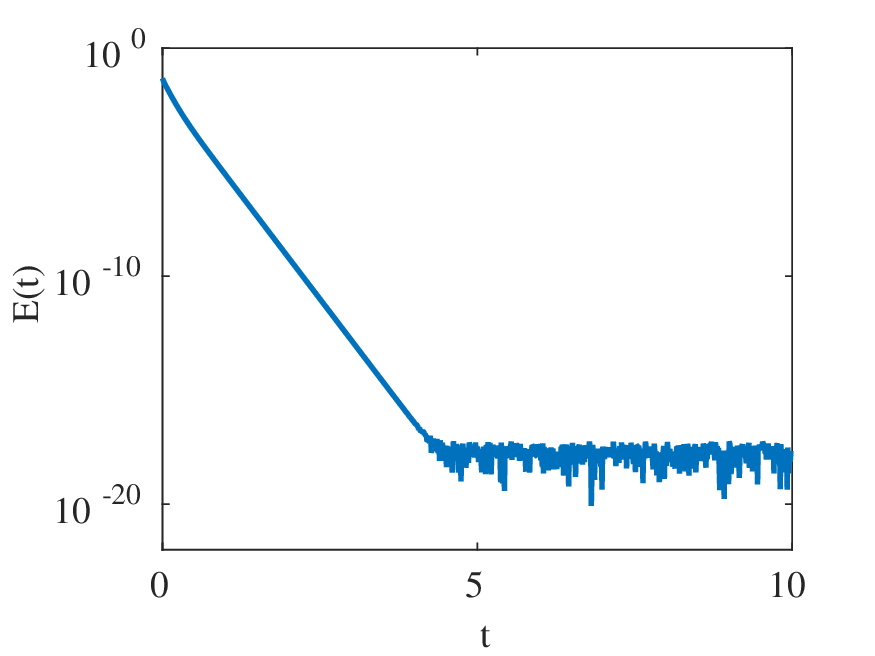}}
\subfigure[$C=0$, $R=R_{SRH}$]{\includegraphics[width=1.6in]{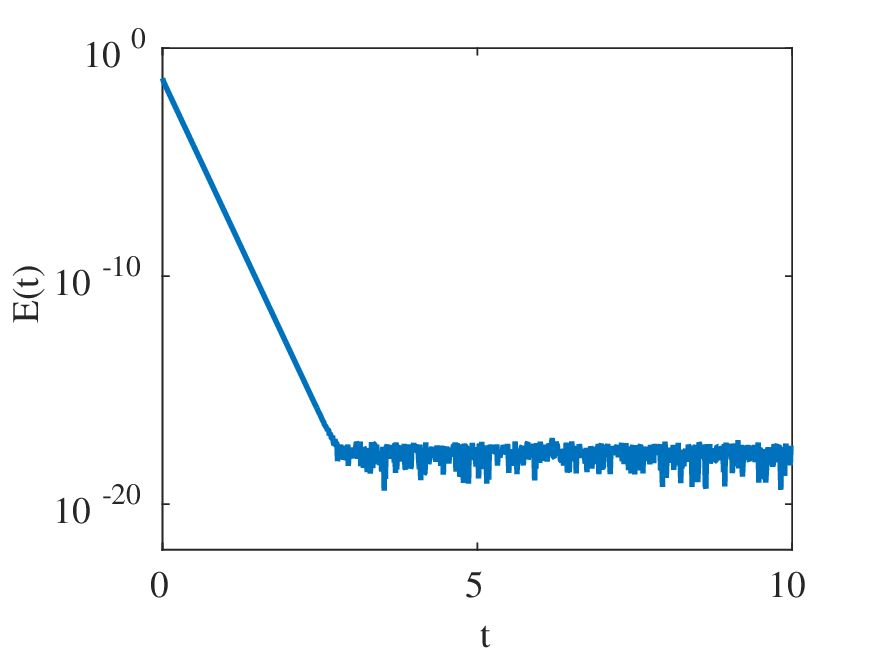}}
\subfigure[$C=0$, $R=R_{AU}$]{\includegraphics[width=1.6in]{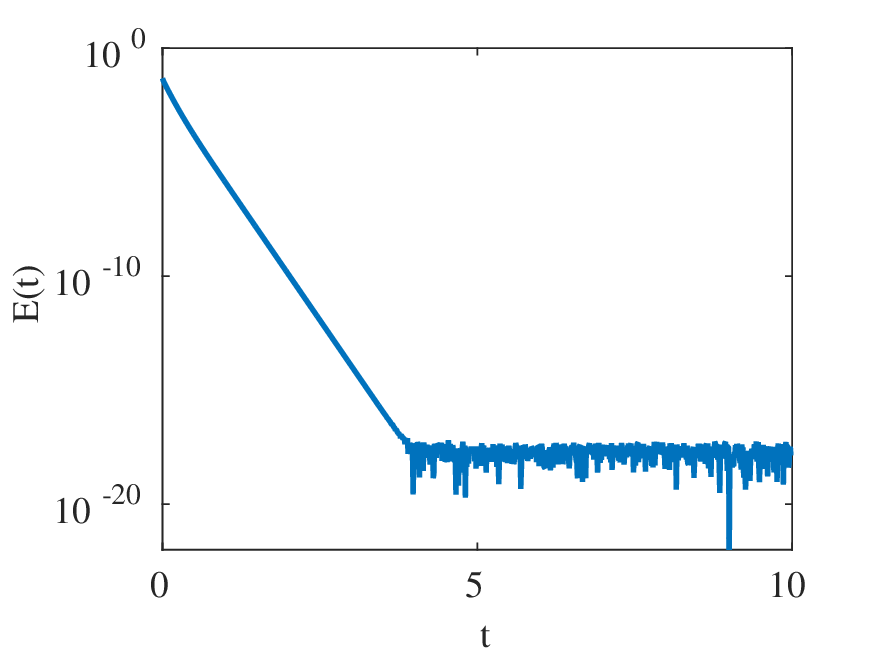}}\\
\subfigure[$C\neq 0$, $R=0$]{\includegraphics[width=1.6in]{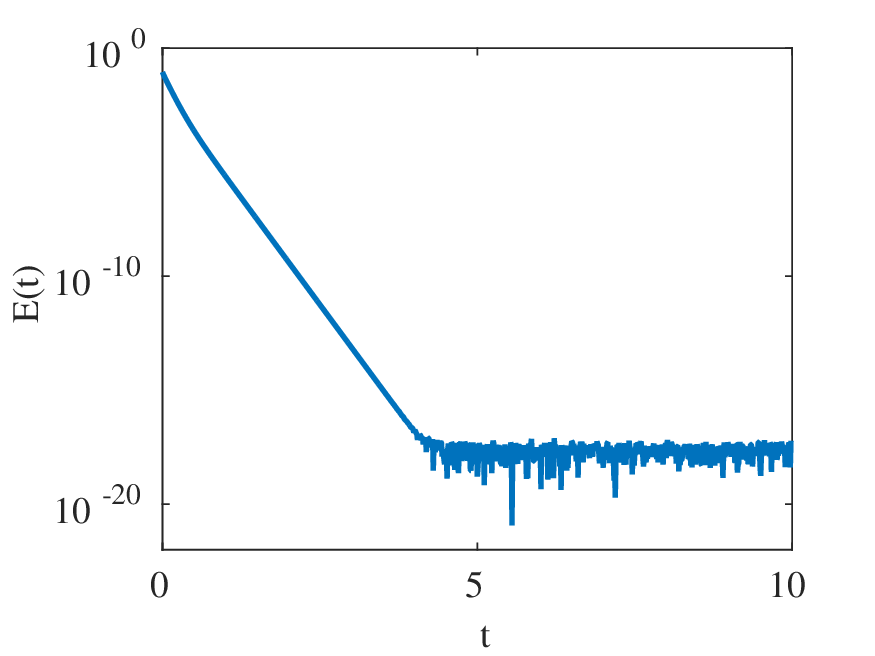}}
\subfigure[$C\neq 0$, $R=R_{SRH}$]{\includegraphics[width=1.6in]{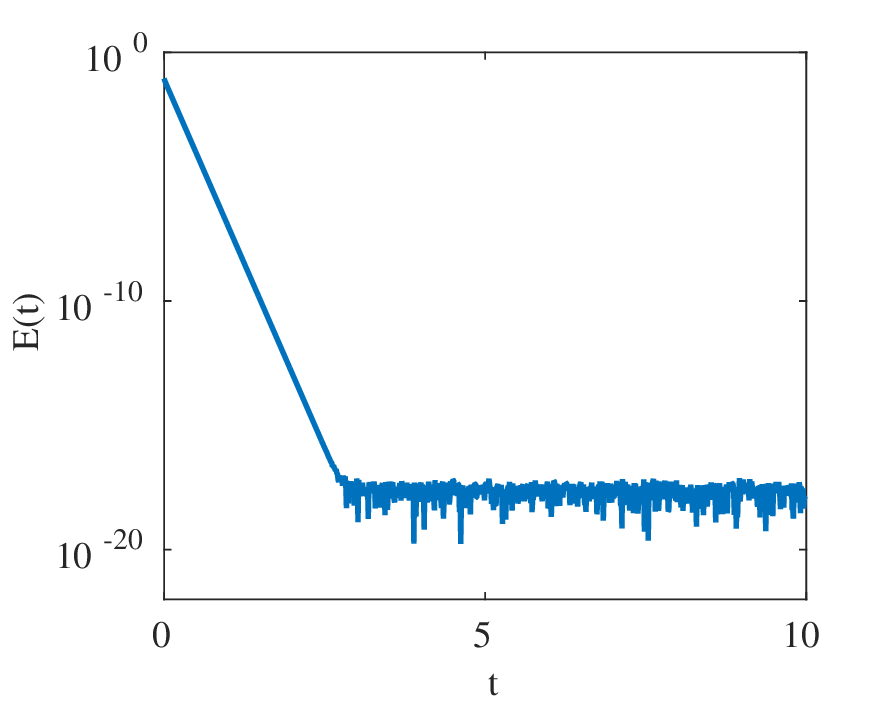}}
\subfigure[$C\neq 0$, $R=R_{AU}$]{\includegraphics[width=1.6in]{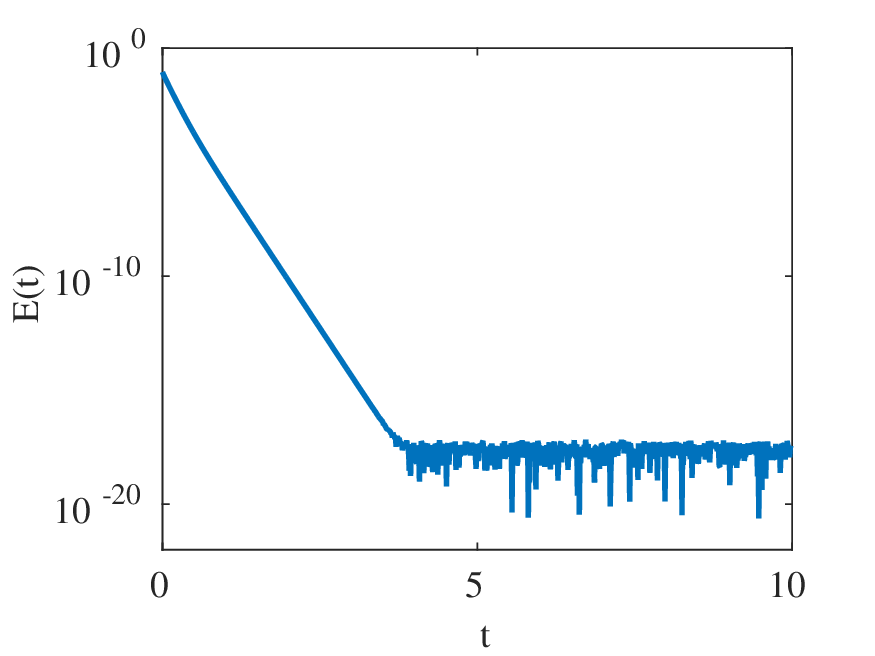}}
\caption{Evolution of $\mathbb{E}^n$ in the linear case with $C=0$ and $C\neq 0$, for different recombination--generation rates.}
\label{fig-lin}
\end{figure}

\begin{figure}
\centering
\subfigure[$C=0$, $R=0$]{\includegraphics[width=1.6in]{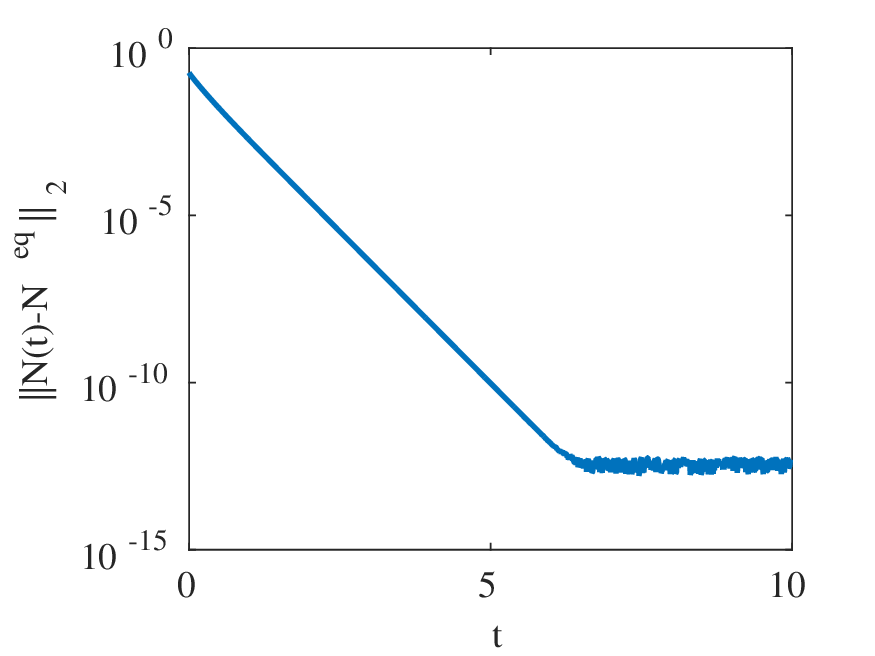}}
\subfigure[$C=0$, $R=R_{SRH}$]{\includegraphics[width=1.6in]{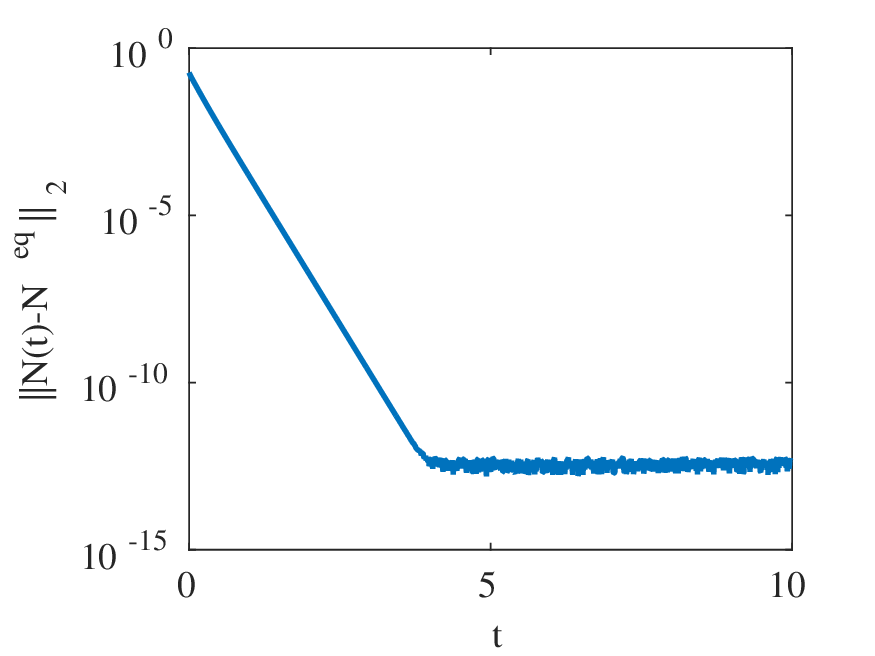}}
\subfigure[$C=0$, $R=R_{AU}$]{\includegraphics[width=1.6in]{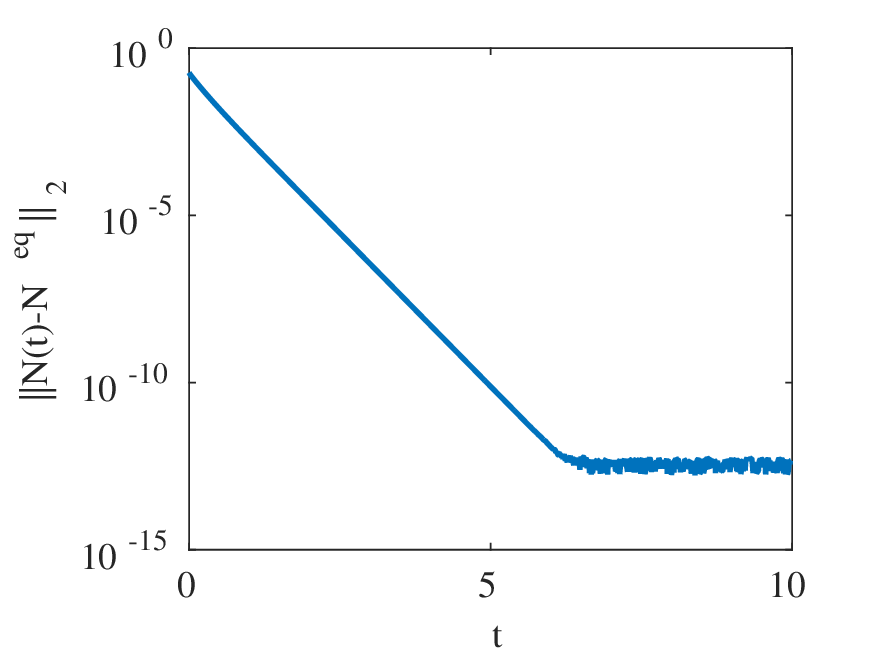}}\\
\subfigure[$C\neq 0$, $R=0$]{\includegraphics[width=1.6in]{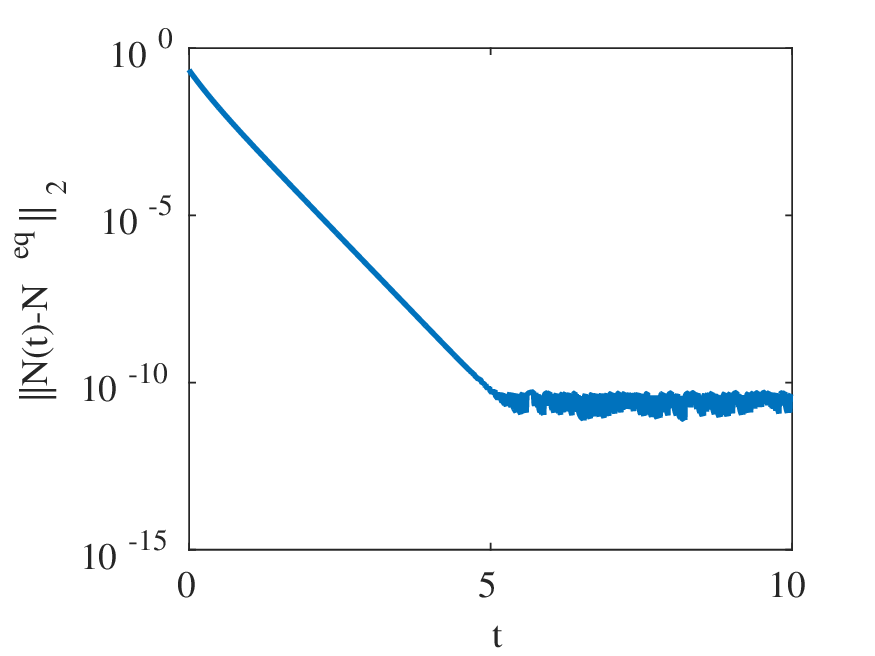}}
\subfigure[$C\neq 0$, $R=R_{SRH}$]{\includegraphics[width=1.6in]{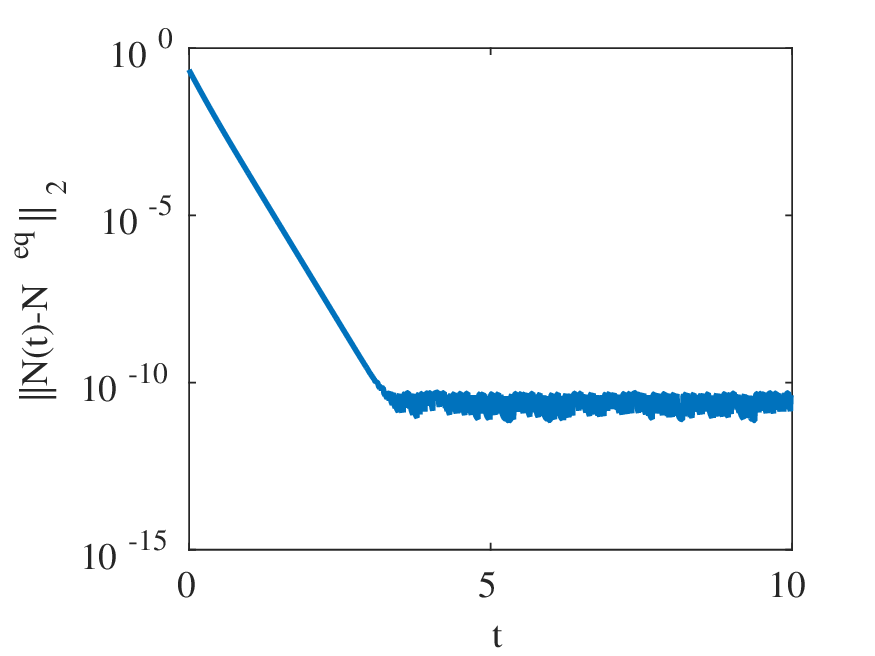}}
\subfigure[$C\neq 0$, $R=R_{AU}$]{\includegraphics[width=1.6in]{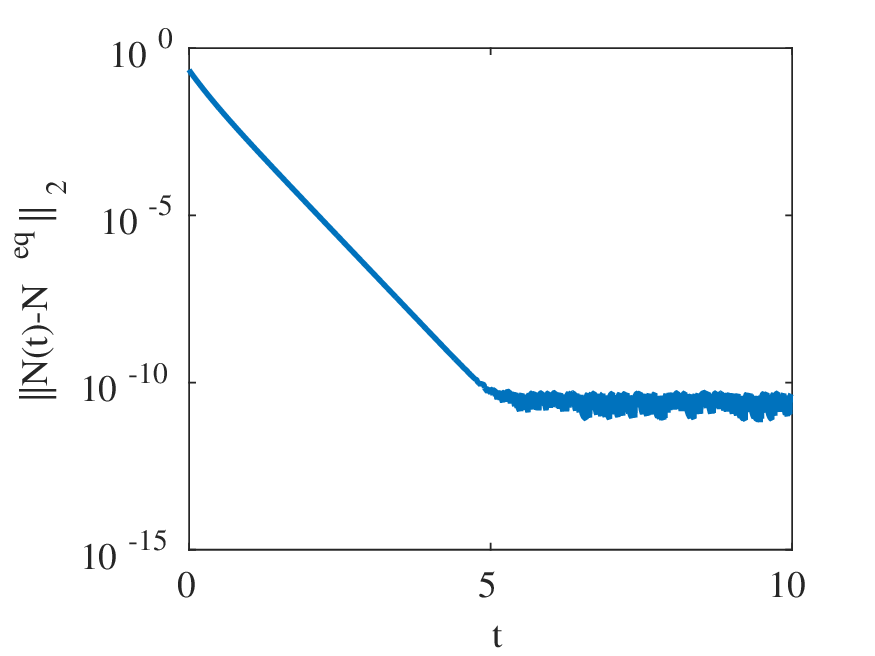}}
\caption{Evolution of $\|N^n-N^{eq}\|_2 $ in the linear case with $C=0$ and $C\neq 0$, for different recombination--generation rates.}
\label{fig-lin-L2}
\end{figure}

\paragraph{Linear case.} We first consider the isothermal case $r=Id$. The boundary conditions are given with $N^D_{0}=e$, $P^D_{0}=e^{-1}$, $N_{1}^D=P_{1}^D=1$. We consider three different recombination--generation rates: 
\begin{itemize}
\item $R=0$,
\item $\ds{R_{SRH}(N,P)=10\,\frac{NP-1}{N+P+1}}$ (Shockley--Read--Hall),
\item $\ds{R_{AU}(N,P)=0.1(N+P)(NP-1)}$ (Auger).
\end{itemize}
In Figure \ref{fig-lin}, we represent the time evolution of the relative entropy in log scale. As expected, we observe an exponential decay rate towards the equilibrium state. The doping does not seem to have an effect on the decay rate, whereas the recombina\-tion--generation rate can modify it. We also represent the evolution of the $L^2$ norm of $N^n-N^{eq}$ in Figure \ref{fig-lin-L2}, which is in good agreement with \eqref{cv-exp-L2-lin}.

\begin{figure}
\centering
\subfigure[$C=0$]{\includegraphics[width=2.4in]{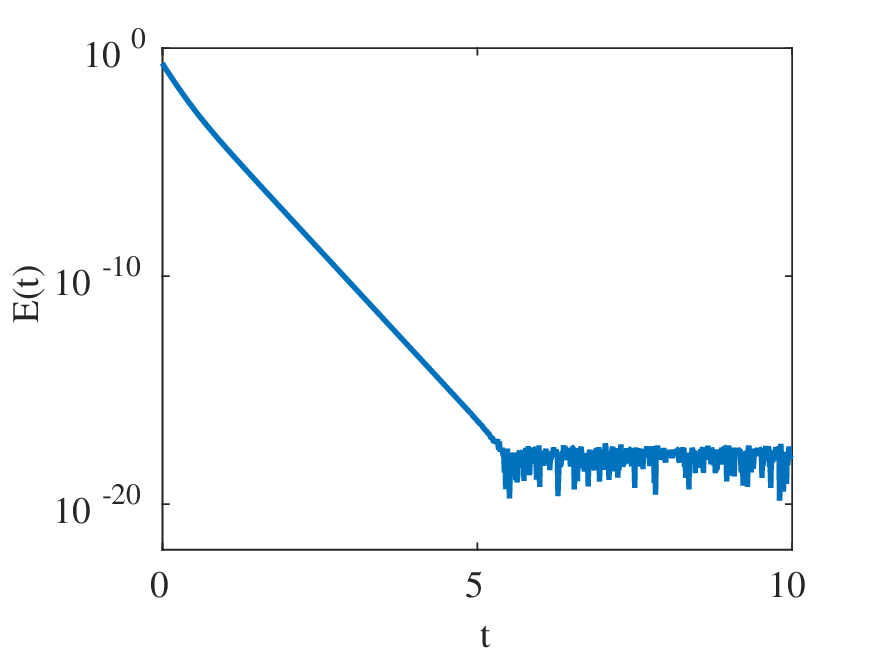}}
\subfigure[$C\neq 0$]{\includegraphics[width=2.4in]{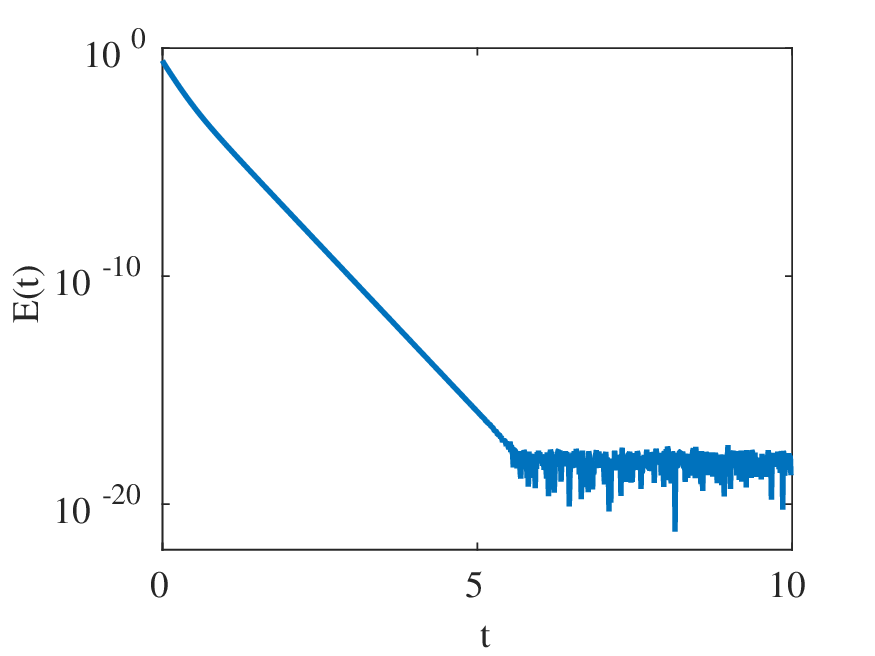}}
\caption{Evolution of $\mathbb{E}^n$ in the nonlinear nondegenerate case, with $C=0$ and $C\neq 0$.}
\label{fig-nonlin-nondeg}
\end{figure}

\begin{figure}
\centering
\subfigure[$C=0$]{\includegraphics[width=2.4in]{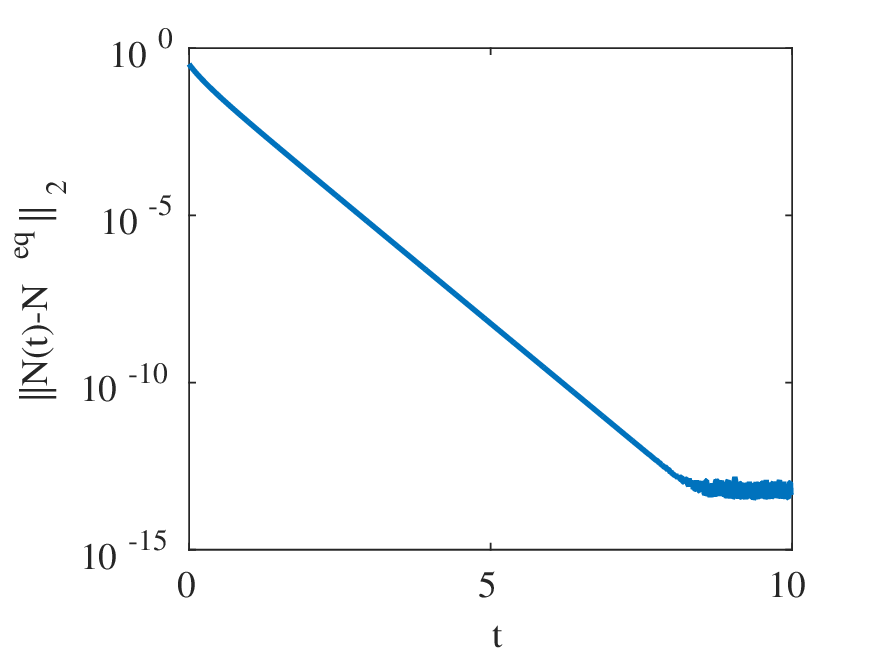}}
\subfigure[$C\neq 0$]{\includegraphics[width=2.4in]{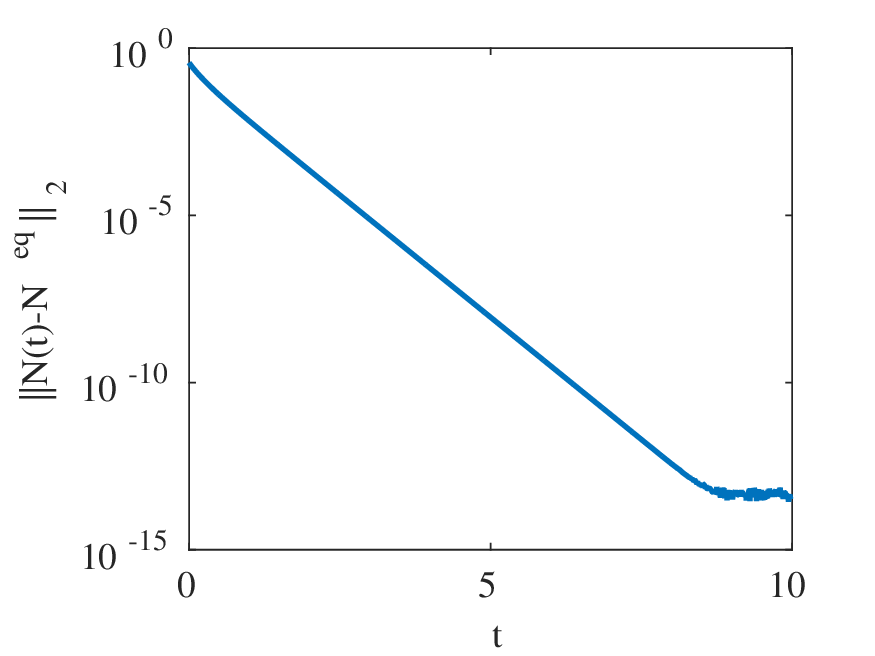}}
\caption{Evolution of $\|N^n-N^{eq}\|_2$ in the nonlinear nondegenerate case, with $C=0$ and $C\neq 0$.}
\label{fig-nonlin-nondeg-L2}
\end{figure}

\begin{figure}
\centering
\includegraphics[width=2.4in]{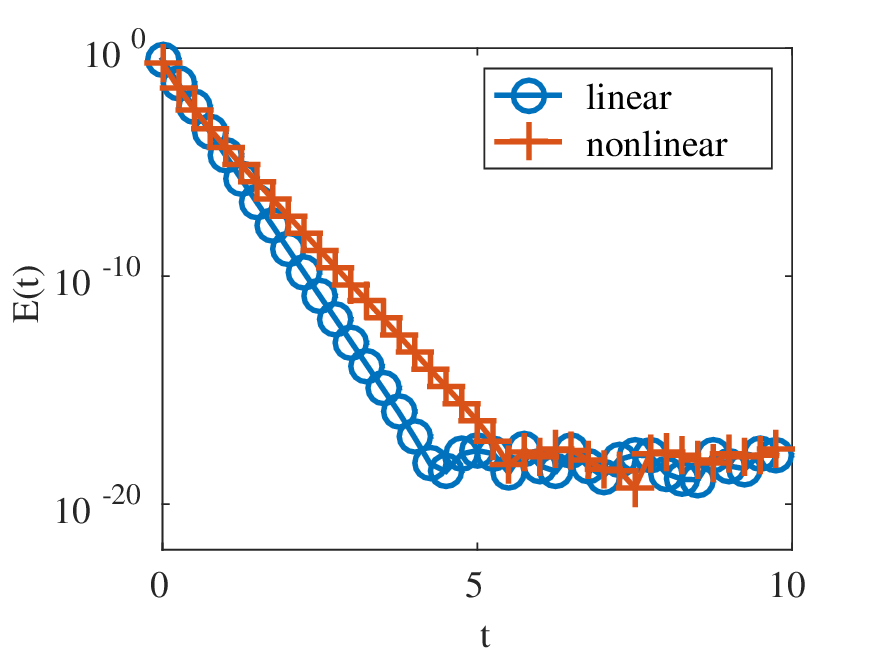}
\caption{Comparison between the linear and the nonlinear cases.}
\label{fig-compare-E}
\end{figure}

\paragraph{Nonlinear nondegenerate case.} We now consider the case of a nonlinear pressure law $r(s)=s^{5/3}$. As mentioned in the introduction, physically relevant recombination--generation rates are not known in this case, then we take $R=0$. The Dirichlet boundary conditions are given with $N^D_{0}=0.9=P^D_{1}$ and $N^D_{1}=0.1=P^D_{0}$. We proved in Theorem \ref{thrm-existence-nonlin} that the approximate densities $N_{\T}^n$, $P_{\T}^n$
satisfy the uniform $L^{\infty}$ estimate \eqref{borne-unif} when $C=0$, and then equations on the densities $N$ and $P$ do not degenerate in this case. In Figures \ref{fig-nonlin-nondeg} and \ref{fig-nonlin-nondeg-L2}, we observe an exponential decay to the thermal equilibrium state, which is in agreement with Theorem \ref{thrm:decayE}. This decay rate is still observed in the case of a nonvanishing doping profile, even it is yet not proved rigorously. In Figure \ref{fig-compare-E}, we compare the relative entropy obtained in this nonlinear case with that obtained with the same data but with $r(s)=s$. The decay rate appears to be slower in the nonlinear case.

\begin{figure}
\centering
\subfigure[$C=0$]{\includegraphics[width=2.4in]{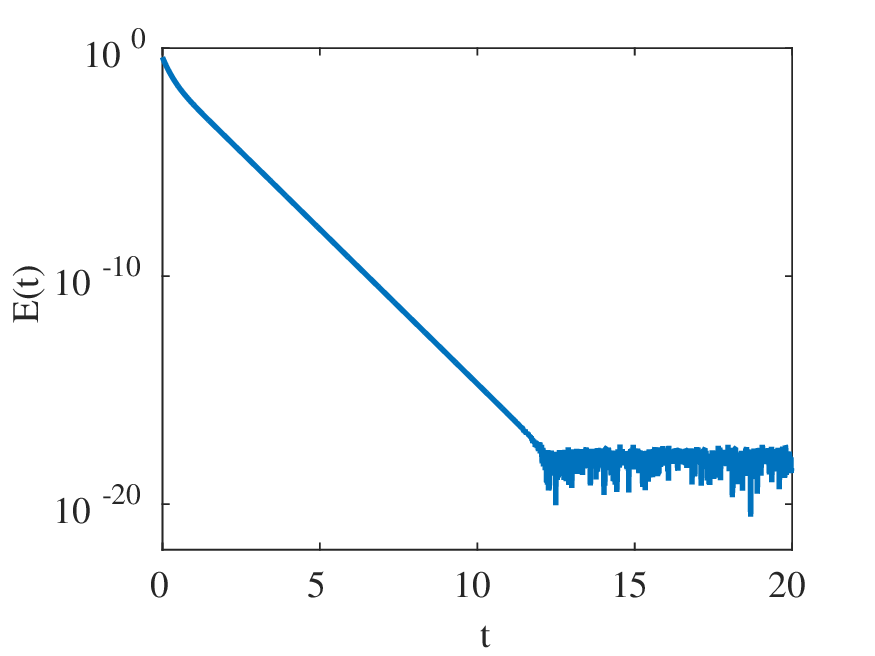}}
\subfigure[$C\neq 0$]{\includegraphics[width=2.4in]{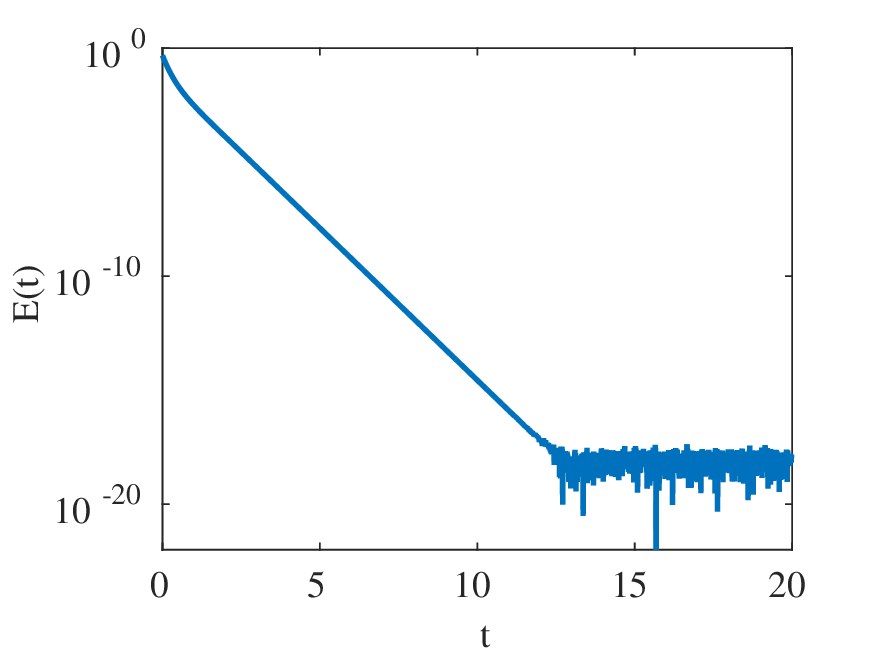}}
\caption{Evolution of $\mathbb{E}^n$ in the nonlinear degenerate case, with $C=0$ and $C\neq 0$.}
\label{fig-nonlin-deg}
\end{figure}

\begin{figure}
\centering
\includegraphics[width=2.4in]{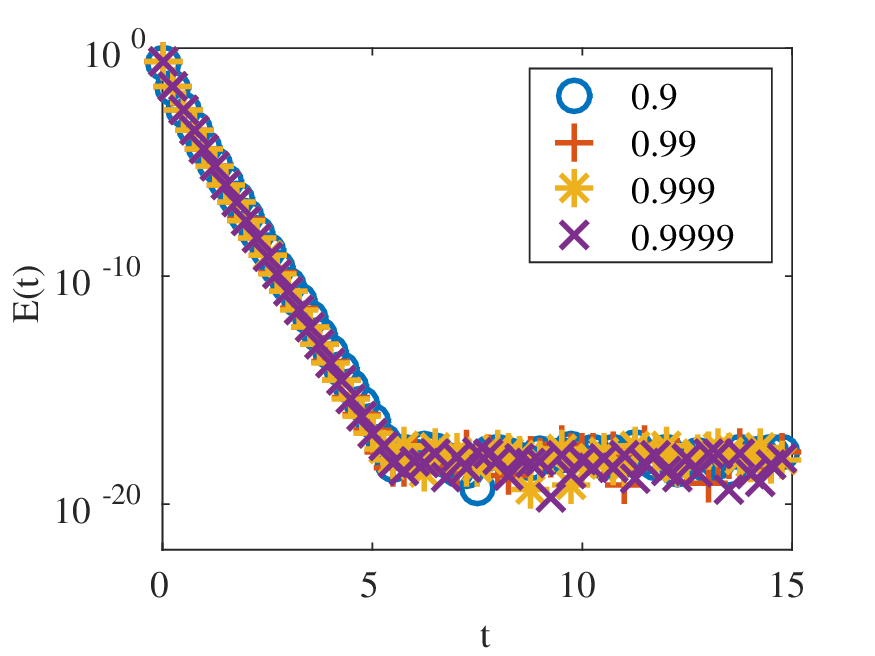}
\caption{Evolution of the relative entropy $\mathbb{E}^n$ with $m=0.1$ and different values of~$M$.}
\label{fig-compare-E2}
\end{figure}

\begin{figure}
\centering
\includegraphics[width=2.4in]{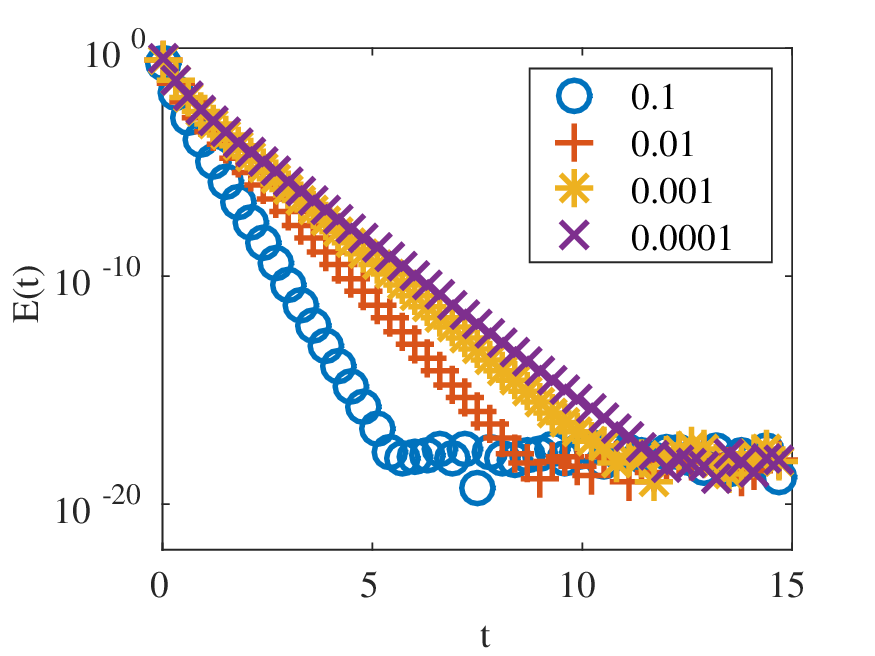}
\caption{Evolution of the relative entropy $\mathbb{E}^n$ with $M=0.9$ and different values of $m$.}
\label{fig-compare-E3}
\end{figure}

\paragraph{Nonlinear degenerate case.} We finally consider the same test case, with Dirichlet conditions vanishing on a part of the boundary: $N^D_{0}=1=P^D_{1}$ and $N^D_{1}=0=P^D_{0}$. We note that the diffusion degenerates when the densities vanish. In this case, we still observe in Figure \ref{fig-nonlin-deg} an exponential convergence of the relative entropy, but with a slower decay rate. In Figure \ref{fig-compare-E2}, we compare the results obtained for $N^D_{1}=0.1=P^D_{0}$ and for different values $M$ of $N^D_{0}=P^D_{1}$: $M=0.9,\,0.99,\,0.999,\,0.9999$. It seems that the variation of this parameter has no influence over the decay rate of the relative entropy. Furthermore, we consider the same test case but with $N^D_{0}=0.9=P^D_{1}$ and different values $m$ of $N^D_{1}=1=P^D_{0}$: $m=0.1,\,10^{-2},\,10^{-3},\,10^{-4}$. In Figure \ref{fig-compare-E3}, it appears that the smaller $m$ is, the slower the decay rate is. This confirms that the slower decay rate observed in Figure \ref{fig-nonlin-deg} is due to the degeneracy of the diffusion.

\begin{figure}
\centering
\subfigure[$C=0$]{\includegraphics[width=2.4in]{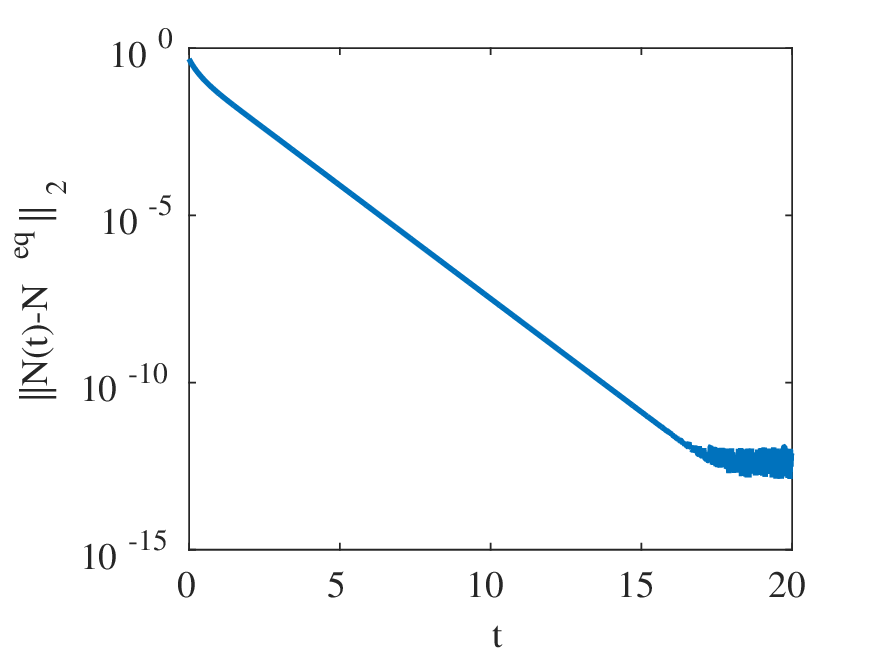}}
\subfigure[$C\neq 0$]{\includegraphics[width=2.4in]{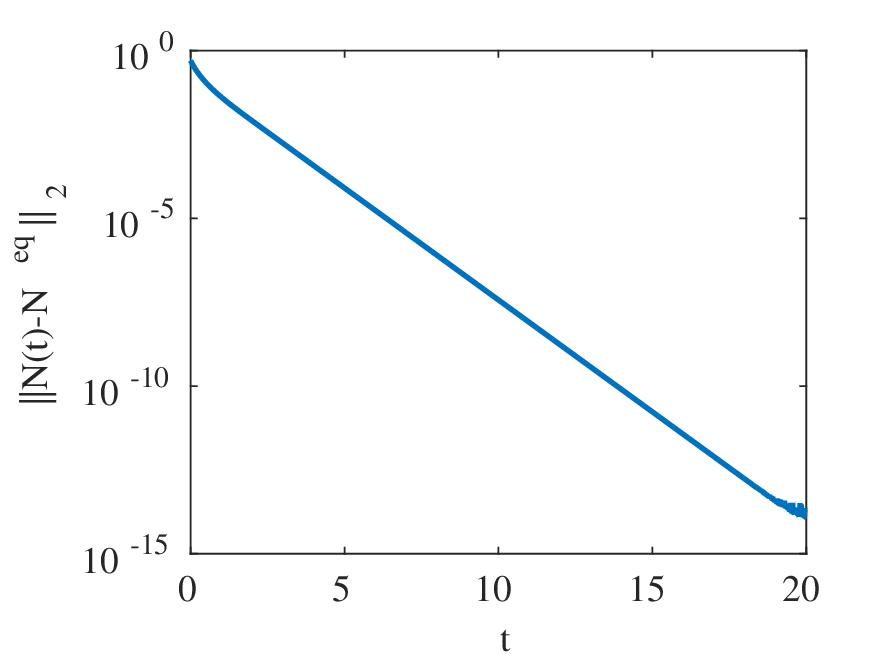}}
\caption{Evolution of $\|N^n-N^{eq}\|_2$ in the nonlinear degenerate case, with $C=0$ and $C\neq 0$.}
\label{fig-nonlin-deg-L2}
\end{figure}

\section{Conclusion}

In this article, we study the large time behavior of a finite volume scheme with Scharfetter--Gummel fluxes discretizing the drift--diffusion model for semiconductors. We prove the convergence of the approximate solution towards an approximation of the thermal equilibrium at an exponential rate as time tends to infinity. This result is established on one hand in the case of a linear diffusion with rather general recombination--generation rate, and on the other hand in the case of a nonlinear diffusion, neglecting recombination and generation processes. In the spirit of \cite{Bessemoulin-Chatard2012}, we consider mostly power functions for the pressure law, corresponding to high density limit of the Fermi-Dirac distribution. Nevertheless, it seems that our result could be applied to more general distribution functions arising in the modeling of organic semiconductors \cite{Koprucki2015}. \\
Moreover, our main theorem is established assuming that uniform-in-time $L^\infty$ estimates hold for the charge carrier densities. This assumption is fulfilled in the case of zero doping profile. Future work would be to prove these uniform-in-time estimates for general $L^\infty$ doping profiles.

\paragraph{Acknowledgements.} The first author thanks the project ANR-12-IS01-0004 GeoNum and the project ANR-14-CE25-0001 Achylles for their partial financial contributions. The second author thanks the team Inria/Rapsodi, the ANR MOONRISE and the Labex CEMPI (ANR-11-LABX-0007-01) for their support.

\bibliographystyle{plain}
\bibliography{bib-DD-exp}

\end{document}